\documentclass[aos,seceqn,nameyear,dvips]{arximspdf}
\usepackage{multirow}
\usepackage{dcolumn}

\doi{10.1214/10-AOS826}
\volume{38}
\issue{6}
\pubyear{2010}
\firstpage{3696}
\lastpage{3723}

\makeatletter

\newcolumntype{k}[1]{D{,}{}{#1}}

\newproclaim{Assumption}{Assumption}
\newproclaim{Study}{Study}

\newcommand{\var}{\operatorname{Var}}
\newcommand{\cov}{\operatorname{Cov}}
\newcommand{\tr}{\operatorname{tr}}
\newcommand{\pr}{\operatorname{Pr}}

\newcommand\independent{\protect\mathpalette{\protect\independenT}{\perp}}
\def\independenT#1#2{\mathrel{\rlap{$#1#2$}\mkern2mu{#1#2}}}

\newtheorem{theorem}{Theorem}
\newtheorem{pro}{Proposition}
\newtheorem{lemma}{Lemma}
\newproclaim{defi}{Definition}

\makeatother

\begin{document}
\begin{frontmatter}

\title{Coordinate-independent sparse sufficient dimension reduction and
variable selection}
\runtitle{Sparse sufficient dimension reduction}

\begin{aug}
\author[A]{\fnms{Xin} \snm{Chen}\thanksref{T1}\ead[label=e1]{ishinchen@gmail.com}},
\author[B]{\fnms{Changliang} \snm{Zou}\thanksref{T2}\ead[label=e2]{chlzou@yahoo.com.cn}} and
\author[C]{\fnms{R. Dennis} \snm{Cook}\corref{}\thanksref{T1}\ead[label=e3]{dennis@stat.umn.edu}}
\runauthor{X. Chen, C. Zou and R. D. Cook}
\affiliation{Syracuse University, Nankai University and University of
Minnesota}
\address[A]{X. Chen\\
 Department of Mathematics\\
Syracuse University\\
Syracuse, New York 13244\\
USA\\
\printead{e1}}
\address[B]{C. Zou\\
LPMC and School\\
\quad of Mathematical Sciences\\
Nankai University\\
Tianjin 300071\\
China\\
\printead{e2}}
\address[C]{R. D. Cook\\
School of Statistics\\
University of Minnesota\\
Minneapolis, Minnesota 55455\\
USA\\
\printead{e3}}
\end{aug}

\thankstext{T1}{Supported by NSF Grants DMS-07-04098 and DMS-10-07547.}
\thankstext{T2}{Supported by NNSF of China Grants 10771107, 11001138 and 11071128.}

\received{\smonth{8} \syear{2009}}
\revised{\smonth{3} \syear{2010}}

%
\begin{abstract}
Sufficient dimension reduction (SDR) in regression, which reduces the
dimension by replacing original predictors with a minimal set of their
linear combinations without loss of information, is very helpful when
the number of predictors is large. The standard SDR methods suffer
because the estimated linear combinations usually consist of all
original predictors, making it difficult to interpret. In this paper,
we propose a unified method---coordinate-independent sparse estimation
(CISE)---that can simultaneously achieve sparse sufficient dimension
reduction and screen out irrelevant and redundant variables
efficiently.
CISE is subspace oriented in the sense that it incorporates a
coordinate-independent penalty term with a broad series of model-based
and model-free SDR approaches. This results in a Grassmann manifold
optimization problem and a fast algorithm is suggested. Under mild
conditions, based on manifold theories and techniques, it can be shown
that CISE would perform asymptotically as well as if the true
irrelevant predictors were known, which is referred to as the oracle
property. Simulation studies and a real-data example demonstrate the
effectiveness and efficiency of the proposed approach.
\end{abstract}

%

\end{frontmatter}

\section{Introduction}

Consider the regression of a univariate response
$y$ on $p$ random predictors ${\mathbf x} = (x_1,\ldots,x_p)^T \in
\mathbb{R}^p$, with the general goal of inferring about
the conditional distribution of $y|{\mathbf x}$. When $p$ is large,
most statistical methods face the ``curse of
dimensionality,'' and thus dimension reduction is desirable.

Sufficient dimension reduction (SDR) introduced by
Cook (\citeyear{Cook1994}, \citeyear{Cook98a}) is important in both
theory and practice. It strives to reduce the
dimension of ${\mathbf x}$ by replacing it with a minimal set of linear
combinations of ${\mathbf x}$, without loss of information on the
conditional distribution of $y|{\mathbf x}$. If a predictor subspace
$\mathcal{S} \subseteq\mathbb{R}^p$ satisfies
\[
y \independent{\mathbf x} | P_{\mathcal{S}}{\mathbf x} ,
\]
where
$\independent$ stands for independence and $P_{(\cdot)}$ represents
the projection matrix with respect to the standard inner product,
then $\mathcal{S}$ is called a dimension reduction space. The
central subspace $\mathcal{S}_{y|{\mathbf x}}$, which is the
intersection of all dimension reduction spaces, is an essential
concept of SDR. Under mild conditions, it can be shown that
$\mathcal{S}_{y|{\mathbf x}}$ is itself a dimension reduction subspace
[Cook (\citeyear{Cook1994}, \citeyear{Cook98a})], which we assume
throughout this article, and
then it is taken as the parameter of interest. The dimension $d$ of
$\mathcal{S}_{y|{\mathbf x}}$, usually far less than $p$, is assumed to
be known in this article. We also assume throughout that $n > p$.

There has been considerable interest in dimension reduction methods
since the introduction of sliced inverse regression [SIR;
\citet{L91}]
and sliced average variance estimation [SAVE; \citet{CW91}].
\citet{L92} and \citet{Cook98b} proposed and studied the method
of principal Hessian directions (PHD), and the related method of
iterative Hessian transformations was proposed by \citet{CL02}.
\citet{CCL02} proposed partial sliced inverse
regression for estimating a partial central subspace. \citet{YC02}
introduced a covariance method for estimating the central
$k$th moment subspace. Most of these and many other dimension
reduction methods are based on the first two conditional moments and
as a class are called F2M methods [\citet{CF09}]. They
provide exhaustive estimation of $\mathcal{S}_{y|{\mathbf x}}$ under
mild conditions. Recently, \citet{LiWang2007} proposed another F2M
method called directional regression (DR). They argued that DR is
more accurate than or competitive with all of the previous F2M
dimension reduction proposals. 
In contrast to these and other
moment-based SDR approaches, \citet{Cook07} introduced a
likelihood-based paradigm for SDR that requires a model for the inverse
regression of ${\mathbf x}$ on $y$. This paradigm, which is broadly
referred to as principal fitted components (PFC), was developed further
by \citet{CF09}. Likelihood-based SDR inherits properties and
methods from
general likelihood theory and can be very efficient in estimating
the central subspace.

All of the aforementioned dimension reduction methods suffer because
the estimated linear reductions usually involve all of the original
predictors ${\mathbf x}$. As a consequence, the results can be hard to
interpret, the important variables may be difficult to identify and
the efficiency gain may be less than that possible with variable
selection. These limitations can be overcome by screening irrelevant
and redundant predictors while still estimating a few linear
combinations of the active predictors. Some attempts have been made
to address this problem in dimension reduction generally and SDR in
particular. For example, \citet{LCN05} proposed a model-free
variable selection method based on SDR. \citet
{ZouHastieTibshirani2006} proposed a
sparse principal component analysis. \citet{NCT05} introduced a
shrinkage version of SIR, while \citet{LN06} suggested a
sparse version of SIR. \citet{Li2007} studied sparse SDR by
adapting the
approach of \citet{ZouHastieTibshirani2006}. \citet
{ZhouHe2008} proposed a
constrained canonical correlation procedure ($C^3$) based on
imposing the $L_1$-norm constraint on the effective dimension
reduction estimates in CANCOR [\citet{Fungetal02}], followed by a
simple variable filtering method. Their procedure is attractive
because they showed that it has the oracle property [\citet{DJ94},
\citet{FanLi2001}]. More recently, \citet{LengWang2009}
proposed a general adaptive sparse principal component
analysis and \citet{JL09} studied the large $p$ theory in
sparse principal components analysis.

However, most existing sparse dimension reduction methods are
conducted stepwise, estimating a sparse solution for a basis matrix
of the central subspace column by column. Instead, in this article,
we propose a unified one-step approach to reduce the number of
variables appearing in the estimate of $\mathcal{S}_{y|{\mathbf x}}$.
Our approach, which hinges operationally on Grassmann manifold
optimization, is able to achieve dimension reduction and variable
selection simultaneously. Additionally, our proposed method has the
oracle property: under mild conditions the proposed estimator would
perform asymptotically as well as if the true irrelevant predictors
were known.\looseness=1

We start in Section~\ref{sec21} by reviewing the link between many SDR
methods and a generalized eigenvalue problem disclosed by \citet{Li2007}.
In Section~\ref{sec22}, we describe a new SDR penalty function that is
invariant under orthogonal transformations and targets the removal
of row vectors from the basis matrix. Based on this penalty
function, in Section~\ref{sec23}, a coordinate-independent penalized
procedure is proposed which enables us to incorporate many
model-free and model-based SDR approaches into a simple and unified
framework to implement variable selection within SDR. A fast
algorithm, which combines a local quadratic approximation
[\citet{FanLi2001}] and an eigensystem analysis in each iteration
step, is
suggested in Section~\ref{sec24} to handle our Grassmann manifold
optimization problem with its nondifferentiable penalty function.
In Section~\ref{sec25}, we describe the oracle property of our estimator. Its
proof differs significantly from those in the context of variable
selection in single-index models [e.g., \citet{FanLi2001},
\citet{Zou2006}]
because the focus here is on subspaces rather than on coordinates.
Results of simulation studies are reported in Section~\ref{sec3}, and the Boston
housing data, is analyzed in Section~\ref{sec4}. Concluding remarks about the
proposed method can be found in Section~\ref{sec5}. Technical details are
given in the \hyperref[app]{Appendix}.\looseness=1

\section{Theory and methodology}\label{sec2}

\subsection{Motivation: Generalized eigenvalue
problems revisited}\label{sec21}

Li (\citeyear{Li2007})\break showed that many moment based sufficient dimension
reduction methods can be formulated as a generalized eigenvalue
problem in the following form:
%
%
\begin{equation}\label{evf}
\mathbf{M}_n \bolds\delta_{ni} = \lambda_{ni} \mathbf{N}_n \bolds
\delta_{ni}\qquad
\mbox{for } i=1,\ldots,p,
\end{equation}
where $\mathbf{M}_n \geq0$ is a method-specific symmetric kernel matrix,
$\mathbf{N}_n > 0$ is symmetric, often taking the form of the sample covariance
matrix $\bolds\Sigma_n$ of ${\mathbf x}$;
$\bolds\delta_{n1}, \ldots,\bolds\delta_{np}$ are eigenvectors
such that $\bolds\delta_{ni}^T \mathbf{N}_n \bolds\delta_{nj} = 1$
if $i = j$
and 0 if $i \neq j$ and $\lambda_{n1} \geq\cdots\geq\lambda_{np}$
are the corresponding eigenvalues. We use the subscript ``$n$'' to
indicate that $\bolds\Sigma_n$, $\mathbf{M}_n$, $\mathbf{N}_n$ and
$\lambda_{ni}$ are the sample versions of the corresponding
population analogs $\bolds\Sigma$, $\mathbf{M}$, $\mathbf{N}$ and
$\lambda_i$. Under certain conditions that are usually imposed only
on the marginal distribution of ${\mathbf x}$, the first $d$
eigenvectors $\{\bolds\delta_{n1}, \ldots,\bolds\delta_{nd}\}$, which
correspond to the nonzero eigenvalues $\lambda_{n1} >\cdots>
\lambda_{nd}$ form a consistent estimator of a basis for the central
subspace. Letting $ \mathbf{z}= {\bolds\Sigma}^{-1/2}\{{\mathbf
x}-E({\mathbf x})\}$. Many commonly used moment based SDR methods are
listed in
Table~\ref{table1} with the population versions of $\mathbf{M}_n$ and~$\mathbf{N}_n$.

%
%
\begin{table}
\caption{The generalized eigenvalue formulations for principle
component analysis (PCA), principle fitted component (PFC) models,
sliced inverse regression (SIR), sliced average variance estimation
(SAVE) and directional regression (DR) methods}\label{table1}
\begin{tabular*}{\tablewidth}{@{\extracolsep{\fill}}lcc@{}}
\hline
\textbf{Method} & \multicolumn{1}{c}{$\mathbf{M}$}  &
\multicolumn{1}{c@{}}{$\mathbf{N}$}\\
\hline
PCA & ${\bolds\Sigma}$ & $\mathbf{I}_p$\\
PFC & ${\bolds\Sigma}_{\mathrm{fit}}$ & ${\bolds\Sigma}$ \\
SIR & $\cov[E\{{\mathbf x}-E({\mathbf x})|y\}]$ & ${\bolds\Sigma}$
\\
SAVE & ${\bolds\Sigma}^{1/2}E[\{\mathbf{I}_p-\cov(\mathbf{z}|y)\}
^2]{\bolds\Sigma
}^{1/2}$ & ${\bolds\Sigma}$ \\
DR & ${\bolds\Sigma}^{1/2}\{2E[E^2(\mathbf{z}\mathbf{z}^T |y)] +
2E^2[E(\mathbf{z}|y)E(\mathbf{z}^T |y)]$ & \\
& \hspace*{24.34pt}\qquad\qquad$+\, 2E[E(\mathbf{z}|y)E(\mathbf{z}|y)]E[E(\mathbf{z}|y)E(\mathbf
{z}^T |y)] -2\mathbf{I}_p\}{\bolds\Sigma}^{1/2}$ &  ${\bolds\Sigma}$
\\ \hline
\end{tabular*}
\end{table}

Following \citet{Cook2004}, \citet{Li2007} showed that the
eigenvectors
$\{\bolds\delta_{n1}, \ldots,\break\bolds\delta_{nd}\}$ from (\ref
{evf}) can be
obtained by minimizing a least square objective function. Let
%
%
\begin{equation}\label{lsobj}
\widehat{\mathbf{V}} =\mathop{\arg\min}_{\mathbf{V}}\sum
_{i=1}^p\|\mathbf{N}_n^{-1}\mathbf{m}_i - \mathbf{V}\mathbf{V}^T
\mathbf{m}_i\|_{\mathbf{N}_n}^2\qquad\mbox{subject to }
\mathbf{V}^T \mathbf{N}_n \mathbf{V} = \mathbf{I}_d,\hspace*{-35pt}
\end{equation}
where $\mathbf{m}_i$ denotes the $i$th column of $\mathbf{M}_n^{1/2}$,
$i=1,\ldots,p$, $\mathbf{V}$ is a $p\times d$ matrix and the norm
here is
with respect to the $\mathbf{N}_n$ inner product. Then
$\widehat{\mathbf V}_j =\bolds\delta_{nj}$, $j=1,\ldots,d$, where
$\widehat{\mathbf V}_j$
stands for the $j$th column of $\widehat{\mathbf V}$, so that
$\mathrm{span}(\widehat{\mathbf V})$ is the estimator of the central
subspace. To get a sparse solution, Li then added penalties to
the objective function in (\ref{lsobj}), leading to the optimization
problem
\[
(\widehat{\bolds\alpha}, \widehat{\mathbf V}_{s}) =
\min_{\bolds\alpha,\mathbf
V}\Biggl\{\sum_{i=1}^p\|\mathbf{N}_n^{-1}\mathbf{m}_i - \bolds\alpha
\mathbf V^T
\mathbf{m}_i\|_{\mathbf{N}_n}^2 + \tau_2 \operatorname{tr}(\mathbf V^T
\mathbf{N}_n \mathbf
V) + \sum_{i=1}^d \tau_{1,j}\|\mathbf V_j\|_1 \Biggr\},
\]
subject to $\bolds\alpha^T \mathbf{N}_n \bolds\alpha= \mathbf{I}_d$,
where $\tr(\cdot)$ stands for the trace operator, \mbox{$\|\cdot\|_r$} denotes
the $L_r$ norm, $\tau_2$ is some positive constant and $\tau_{1,j}
\geq0$ for $j=1,\ldots,d$ are the lasso shrinkage parameters that
need to be determined by some method like cross validation (CV). The
solution $\widehat{\mathbf V}_{s}$ is called the sparse sufficient
dimension reduction estimator. As a result of the lasso constraint,
$\widehat{\mathbf V}_{s}$ is expected to have some elements shrunk to
zero.

We can see that Li's sparsity method is coordinate dependent
because the $L_1$ penalty term is not invariant under the orthogonal
transformation of the basis and it forces individual elements of the
basis matrix $\widehat{\mathbf V}_{s}$ to zero. However, variable
screening requires that entire rows of $\widehat{\mathbf V}_{s}$ be
zero, which is not the explicit goal of Li's method.
To see this more
clearly, partition ${\mathbf x}$ as $({\mathbf x}_1^T,{\mathbf
x}_2^T)^T$, where
${\mathbf x}_1$ corresponds to $q$ elements of ${\mathbf x}$ and
${\mathbf x}_2$
to the remaining elements. If
%
%
\begin{equation}\label{conind}
y \independent{\mathbf x}_2 | {\mathbf x}_1 ,
\end{equation}
then ${\mathbf x}_2$ can be removed, as given ${\mathbf x}_1$,
${\mathbf x}_2$
contains no further information about~$y$. Let the $p \times d$
matrix $\bolds\eta$ be a basis for $\mathcal{S}_{y|{\mathbf x}}$ and
partition $\bolds\eta= (\bolds\eta_1^T,\bolds\eta_2^T)^T$ in
accordance with
the partition of ${\mathbf x}$. Then the condition (\ref{conind}) is
equivalent to $\bolds\eta_2 = 0$ [\citet{Cook2004}], so the
corresponding rows
of the basis are zero vector.

In effect, Li's method is designed for element screening, not
variable screening. Our experience reflects this limitation and
reinforces the notion that $\widehat{\mathbf V}_{s}$ may not be
sufficiently effective at variable screening. Inspired by Li's
method, we propose a new variable screening method---called
coordinate-independent sparse estimation (CISE)---in the next
subsection. We will show that CISE is simpler and more effective
than Li's method at variable screening.

CISE can be applied not only to moment based SDR approaches but also
model based approaches. \citet{Cook07} and \citet{CF08}
developed several powerful model-based dimension reduction
approaches, collectively referred to as principal fitted components
(PFC). PFC-based SDR methods can also be formulated in the same way
as (\ref{evf}), as summarized in the next proposition. In
preparation, consider the following model for the conditional
distribution of ${\mathbf x}$ given $y$,
%
%
\begin{equation}\label{pfcmodel}
{\mathbf x}=\bolds\mu+ \bolds\Gamma\bolds\xi\mathbf{f}(y) +
\bolds\Delta^{1/2}
\bolds\epsilon,
\end{equation}
where $\bolds\mu\in{\mathbb R}^{p}$ is a location vector, $\bolds
\Gamma\in
\mathbb{R}^{p \times d}$, $\bolds\Gamma^T\bolds\Gamma=\mathbf
{I}_d$, $\bolds\xi
\in\mathbb{R}^{d \times r}$ with rank $d$, $\mathbf{f} \in
\mathbb{R}^r$ is a known vector-valued function of $y$, $\bolds\Delta=
\operatorname{Var}({\mathbf x}|y) > 0$, and $\bolds\epsilon\in\mathbb
{R}^p$ is
assumed to be independent of $y$ and normally distributed with mean
0 and identity covariance matrix.
\begin{pro}\label{prop1}
Suppose the conditional distribution of ${\mathbf x}$ given $y$ can be
described by (\ref{pfcmodel}). Then the maximum likelihood estimator
(MLE) of $\mathcal{S}_{y|{\mathbf x}}$ can be obtained through the
generalized eigenvalue problem of the form (\ref{evf}) with
$\mathbf{M}_n =\widehat{\bolds\Sigma}_{\mathrm{fit}}$ and
$\mathbf{N}_n =\bolds\Sigma_n$, where
$\widehat{\bolds\Sigma}_{\mathrm{fit}}$ is the sample covariance matrix
of the fitted vectors from the linear regression of ${\mathbf x}$ on
$\mathbf{f}$.
\end{pro}

A commonly\vspace*{1pt} used case in the PFC models is $\bolds\Delta=\sigma
^{2}\mathbf{I}_p$ for $\sigma>0$, in which the MLE of $\mathcal
{S}_{y|{\mathbf x}}$
can be obtained through (\ref{evf}) with $\mathbf{M}_n
=\widehat{\bolds\Sigma}_{\mathrm{fit}}$ and $\mathbf{N}_n =\mathbf{I}_p$.
The covariates $\mathbf{f}(y)$ in model (\ref{pfcmodel}) usually take
form of polynomial, piecewise or Fourier basis functions. Thus, the
PFC models can effectively deal with the nonlinear relationship
between the predictors and the response.

\subsection{A coordinate-independent penalty function}\label{sec22}

Let $\mathbf{V} = (\mathbf{v}_1,\ldots,\mathbf{v}_p)^T$ denote a $p
\times
d$ matrix with rows $\mathbf{v}_i^T$, $i=1,\ldots,p$. In this section,
we introduce a coordinate-independent penalty, depending only on the
subspace spanned by the columns of $\mathbf V$. Let $\mathbf{q}_i$ be the
vector in $\mathbb{R}^p$ with the $i$th component one, else zero.


We define a general coordinate-independent penalty function as
\[
\phi(\mathbf{V}) = \sum_i \theta_i h_i({\mathbf{q}_i^T \mathbf
{V}\mathbf{V}^T\mathbf{q}_i}),
\]
where $\theta_i \geq0$ serve as penalty parameters, and $h_i$ are
positive convex functions defined in $\mathbb{R}^d$. To achieve
variable screening, the functions $h_i$ must be nondifferentiable
at the zero vector. It is clear that the function $\phi$ is
independent of the basis used to represent the span of $\mathbf{V}$,
since for any orthogonal matrix $\mathbf O$, $\phi(\mathbf{V}) =
\phi(\mathbf{V O})$. In fact, any penalty function defined on
$\mathbf{V}\mathbf{V}^T$
meets our requirement.

Given $h_1 =\cdots= h_p = \sqrt{(\cdot)}$, we have a special
coordinate-independent penalty function:
%
%
\begin{equation}\label{penalty}
\rho(\mathbf{V}) = \sum_{i=1}^p \theta_i \|\mathbf{v}_i\|_2.
\end{equation}
A method for selecting the tuning parameters will be discussed in
Section~\ref{sec26}. We can see that the penalty function $\rho$ has the
same form as the group lasso proposed by \citet{YuanLin2006} but
their concepts and usages are essentially different. Through this
article, we shall use only $\rho$ in application and theory to
demonstrate our ideas.

Penalty (\ref{penalty}) is appealing for variable selection because
it is independent of the basis used to represent the span of $\mathbf
V$, $\rho(\mathbf{V}) = \rho(\mathbf{VO})$ for any orthogonal
matrix~$\mathbf O$, and because it groups the row vector coefficients of
$\mathbf V$.
This motivated us to consider the regularized function
(\ref{penalty}) that can shrink the corresponding row vectors of
irrelevant variables to zero. Another appealing feature of using
this penalty is its oracle property, which is discussed in Section
\ref{sec25}.

\subsection{Coordinate-independent sparse estimation}\label{sec23}
Recall the generalized eigenvalue problem (\ref{evf}) and the associated
notation. Formally,
\[
\sum_{i=1}^p\Vert
\mathbf{N}_n^{-1}\mathbf{m}_i-\mathbf{V}\mathbf{V}^T\mathbf
{m}_i\Vert^2_{\mathbf{N}_n}
=\tr(\mathbf{G}_n)-\tr(\mathbf{V}^T\mathbf{M}_n\mathbf{V}),
\]
where $\mathbf{G}_n=\mathbf{N}_n^{-1/2}\mathbf{M}_n\mathbf
{N}_n^{-1/2}$ and we
use $\mathbf{G}$ to denote its population analog in what follows.
Hence, the ordinary sufficient dimension reduction estimation
(OSDRE) given in (\ref{lsobj}) is
%
%
\begin{equation}\label{Vhat}
\widehat{\mathbf{V}}=\mathop{\arg\min}_{\mathbf{V}}-\tr(\mathbf
{V}^T\mathbf{M}_n\mathbf{V})\qquad\mbox{subject to } \mathbf{V}^T\mathbf
{N}_n\mathbf{V}=\mathbf{I}_d.
\end{equation}
%
By using the coordinate independent penalty function given in last
subsection, we propose the following coordinate-independent sparse
sufficient dimension reduction estimator (CISE):
%
%
\begin{equation} \label{Vtilde}
\tilde{\mathbf{V}}=\mathop{\arg\min}_{\mathbf{V}}\{-\tr(\mathbf
{V}^T\mathbf{M}_n\mathbf{V})+\rho(\mathbf{V})\}\qquad\mbox{subject to }
\mathbf{V}^T\mathbf{N}_n\mathbf{V}=\mathbf{I}_d,\hspace*{-35pt}
\end{equation}
where $\rho(\mathbf{V})$ is defined in (\ref{penalty}).

The solution $\tilde{\mathbf{V}}$ is not unique as
$\tilde{\mathbf{V}}\mathbf
O$ is also a solution for any orthogonal matrix~$\mathbf O$. In a strict
sense, we are minimizing (\ref{Vtilde}) over the span of the columns
of~$\mathbf{V}$. Thus, $\tilde{\mathbf{V}}$ denotes any basis of the solution
of (\ref{Vtilde}). Analogously, the solution $\widehat{\mathbf{V}}$ is
one basis of the solution of (\ref{Vhat}). Before proceeding, we
rewrite (\ref{Vhat}) and (\ref{Vtilde}) as equivalent unitary
constrained optimization problems which will facilitate our
exposition. We summarize the result into the following proposition
without giving its proof since it follows from some straightforward
algebra.
\begin{pro}\label{prop2}
The minimizer (\ref{Vhat}) is equivalent to $\widehat{\mathbf
{V}}=\mathbf{N}_n^{-1/2}\widehat{\bolds\Gamma}$ where
%
%
\begin{equation}\label{GammaHat}
\widehat{\bolds\Gamma}=\mathop{\arg\min}_{\bolds\Gamma}-\tr
(\bolds\Gamma^T\mathbf{G}_n\bolds\Gamma)\qquad \mbox{subject to }
\bolds\Gamma^T\bolds\Gamma=\mathbf{I}_d.
\end{equation}
Furthermore, $\mathbf{G}_n\widehat{\bolds\Gamma}=\widehat{\bolds
\Gamma}{\bolds\Lambda}_{n1}$, where
${\bolds\Lambda}_{n1} = \operatorname{diag}({\lambda}_{n1}, \ldots
,{\lambda}_{nd})$.
Correspondingly, the minimizer (\ref{Vtilde}) is
equivalent to $\tilde{\mathbf{V}}=\mathbf{N}_n^{-1/2}\tilde{\bolds
\Gamma}$,
where
%
%
\begin{equation}\label{GammaTilde}\qquad
\tilde{\bolds\Gamma}=\mathop{\arg\min}_{\bolds\Gamma}\{-\tr
(\bolds\Gamma^T\mathbf{G}_n\bolds\Gamma)+\rho(\mathbf
{N}_n^{-1/2}\bolds\Gamma)\}\qquad\mbox{subject
to } \bolds\Gamma^T\bolds\Gamma=\mathbf{I}_d.
\end{equation}
\end{pro}

The minimization of (\ref{GammaHat}) and (\ref{GammaTilde}) is a
Grassmann manifold optimization problem. A Grassmann manifold, which
is defined as the set of all $d$-dimensional subspaces in
$\mathbb{R}^p$, is the natural parameter space for the $\bolds\Gamma$
parametrization in (\ref{GammaHat}). For more background on
Grassmann manifold optimization, see \citet{EAS98}. The
traditional Grassmann manifold optimization techniques cannot be
applied directly to (\ref{GammaTilde}) due to the
nondifferentiability of $\rho(\cdot)$. Nevertheless, we have
devised a simple and fast algorithm to solve (\ref{GammaTilde}), as
discussed in the next subsection.

\subsection{Algorithm}\label{sec24}

To overcome the nondifferentiability of $\rho(\cdot)$, we adopt the
local quadratic approximation of \citet{FanLi2001}; that is, we
approximate the penalty function locally with a quadratic function
at every step of the iteration as follows.

Let $\tilde{\mathbf{V}}^{(0)} = (\tilde{\mathbf{v}}_1^{(0)},\ldots
,\tilde{\mathbf{v}}_p^{(0)})^T = \mathbf{N}_n^{-1/2}\tilde{\bolds
{\Gamma}}{}^{(0)}$ be the starting value.
The unconstrained first derivative of $\rho(\mathbf{V})$ with respect to
the $p\times d$ matrix $\mathbf{V}$ is given by
\[
\frac{\partial\rho}{\partial\mathbf{V}} =
\operatorname{diag}\biggl(\frac{\theta_1}{\|\mathbf{v}_1\|_2},\ldots,\frac
{\theta_i} {\|\mathbf{v}_i\|_2},\ldots,\frac{\theta_p}{\|\mathbf
{v}_p\|_2}\biggr)\mathbf{V}.
\]

Following Fan and Li, the first derivative of $\rho(\mathbf{V})$ around
$\tilde{\mathbf{V}}^{(0)}$ can be approximated by
\[
\frac{\partial\rho}{\partial\mathbf{V}} \approx
\operatorname{diag}\biggl(\frac{\theta_1}{\|\tilde{\mathbf{v}}_1^{(0)}\|
_2},\ldots,\frac{\theta_i} {\|\tilde{\mathbf{v}}_i^{(0)}\|
_2},\ldots,\frac{\theta_p}{\|\tilde{\mathbf{v}}_p^{(0)}\|
_2}\biggr)\mathbf{V}
:=\mathbf{H}^{(0)}\mathbf{V}.
\]
By using the second-order Taylor expansion and some algebraic
manipulation, we have
\[
\rho(\mathbf{V}) \approx\tfrac{1}{2} \tr\bigl(\mathbf{V}^T \mathbf
{H}^{(0)}\mathbf{V}\bigr) + C_0 = \tfrac{1}{2}\tr\bigl( \bolds{\Gamma}^T
\mathbf{N}_n^{-1/2} \mathbf{H}^{(0)} \mathbf{N}_n^{-1/2}
\bolds{\Gamma}\bigr) + C_0,
\]
where $C_0$ stands for a constant with respect to $\mathbf{V}$.

Then find $\tilde{\bolds{\Gamma}}{}^{(1)}$ by minimizing:
\begin{eqnarray*}
&&-\tr(\bolds\Gamma^T\mathbf{G}_n\bolds\Gamma)+ \tfrac{1}{2}\tr
\bigl(\bolds{\Gamma}^T
\mathbf{N}_n^{-1/2} \mathbf{H}^{(0)} \mathbf{N}_n^{-1/2}
\bolds{\Gamma}\bigr)\\
&&\qquad
= \tr\bigl\{\bolds\Gamma^T \bigl(-\mathbf{G}_n+ \tfrac{1}{2}\mathbf
{N}_n^{-1/2} \mathbf{H}^{(0)} \mathbf{N}_n^{-1/2}\bigr) \bolds
{\Gamma}\bigr\}.
\end{eqnarray*}
This minimization problem can be easily solved by the eigensystem
analysis of the matrix $\mathbf{G}_n - 2^{-1}\mathbf{N}_n^{-1/2}
\mathbf{H}^{(0)} \mathbf{N}_n^{-1/2}$, that is, the columns of
$\tilde{\bolds{\Gamma}}{}^{(1)}$ are the first $d$ principal component
directions of $\mathbf{G}_n-2^{-1}\mathbf{N}_n^{-1/2} \mathbf
{H}^{(0)} \mathbf{N}_n^{-1/2}$. Next, let $\tilde{\mathbf
{V}}{}^{(1)}=\mathbf{N}_n^{-1/2}\tilde{\bolds{\Gamma}}{}^{(1)}$ and
start the second round of
approximation of $\rho(\mathbf{V})$. The procedures repeat until it
converges. During the iterations, if $\|\tilde{\mathbf{v}}_i^{(k)}\|_2
\approx0$, say $\|\tilde{\mathbf{v}}_i^{(k)}\|_2<\epsilon$ where
$\epsilon$ is a prespecified small positive number (e.g.,
$\epsilon=10^{-6}$), then the variable $x_i$ is removed.

With respect to the choice of the initial values
$\tilde{\bolds{\Gamma}}{}^{(0)}$, a simple but effective solution is to
use $\tilde{\bolds{\Gamma}}{}^{(0)}=\widehat{\bolds{\Gamma}}$, the minimizer
of (\ref{GammaHat}). With $\widehat{\bolds{\Gamma}}$ as the initial
values, we found that the frequency of nonconvergence is negligible
in all of our simulation studies and the convergence is quite fast,
usually requiring a few dozen iterations. A Matlab interface was
used to implement this CISE algorithm. The programs can be obtained
from the first author upon request.

\subsection{Oracle property}\label{sec25} In what follows, without loss of
generality, we assume
that only the first $q$ predictors are relevant to the regression,
where $d \leq q < p$. Given a $p\times d$ matrix $\mathbf{K}$,
$\mathbf{K}_{(q)}$ and $\mathbf{K}_{(p-q)}$ indicate the sub-matrices
consisting
of its first $q$ and remaining $p-q$ rows. If $\mathbf{K}$ is $p\times
p$, then the notation indicates its first $q$ and the last $p-q$
block sub-matrices. In the context of the single-index model,
\citet{FanLi2001} and \citet{Zou2006} have shown that, with
the proper choice
of the penalty functions and regularization parameters, the
penalized likelihood estimators have the oracle property. With
continuous penalty functions, the coefficient estimates that
correspond to insignificant predictors must shrink toward 0 as the
penalty parameter increases, and these estimates will be exactly 0
if that parameter is sufficiently large. In this section, we present
theorems which establish the oracle property of CISE.

Let $a_n=\max\{\theta_j,j\leq q\}$ and $b_n=\min\{\theta_j,j>q\}$,
where the $\theta_j$'s are the penalty parameters defined in Section
\ref{sec22}, let $\lambda_1\geq\cdots\geq\lambda_p\geq0$ denote the
eigenvalues of $\mathbf{G}$, and define the matrix norm $\Vert\mathbf
{V}\Vert_s=\sqrt{\tr(\mathbf{V}^T \mathbf{V})}$. We also require a
metric $D$
in the set of all subspaces of $\mathbb{R}^p$.
\begin{defi}\label{defi1}
The distance between the subspaces spanned by the columns of $\mathbf
{V}_n$ and $\mathbf{V}$, denoted as $D(\mathbf{V}_n,\mathbf{V})$, is defined
as the square root of the largest eigenvalue of
\[
(P_{\mathbf{V}_n} - P_{\mathbf{V}})^T (P_{\mathbf{V}_n} - P_{\mathbf{V}}).
\]
\end{defi}

This distance criterion was first used by \citet
{LiZhaChiaromonte2005} in the
sufficient dimension reduction setting. See \citet
{GohbergLancasterRodman2006}
for more details. We use the following assumptions to establish the
oracle property.
\begin{Assumption}\label{Assum1}
Let $\mathbf{V}_0$ denote the minimizer of (\ref{Vhat})
when the population matrices $\mathbf{M}$ and $\mathbf{N}$ are used
in place of $\mathbf{M}_n$ and $\mathbf{N}_n$. Then ${\mathbf
{V}_0}_{(p-q)} =
0$.
\end{Assumption}
\begin{Assumption}\label{Assum2}
$\mathbf{M}_n=\mathbf{M}+O_p(n^{-1/2})$ and $\mathbf{N}_n=\mathbf
{N}+O_p(n^{-1/2})$.
\end{Assumption}

Given some mild method-specified conditions, the minimizer of
(\ref{Vhat}) $\widehat{\mathbf{V}}$ is a consistent estimator of a
basis of the central subspace. For example, SIR provides the
consistent estimate of the central subspace given that the linearity
and coverage conditions hold [\citet{Cook98a}, \citet{CCL02}].
Consequently, the population version $\mathbf{V}_0$ will be a basis of
the central subspace. Therefore, Assumption~\ref{Assum1} is a
reasonable one
which facilitates our following presentations. Assumption~\ref{Assum2}
is mild
and typically holds. These two assumptions suffice for our main
results.

We state our theorems here, but their proofs are relegated to the
\hyperref[app]{Appendix}. The constrained objective function in the minimization
problem (\ref{Vtilde}) is denoted as $Q(\mathbf{V};\mathbf{M}_n):=
f(\mathbf{V};\mathbf{M}_n)+\rho(\mathbf{V})$ where $f(\mathbf
{V};\mathbf{M}_n)=-\tr(\mathbf{V}^T\mathbf{M}_n\mathbf{V})$. The
first theorem establishes existence of
CISE.
\begin{theorem}\label{theo1}
If Assumptions~\ref{Assum1} and~\ref{Assum2} hold, $\lambda
_d>\lambda_{d+1}$ and
$\sqrt{n}a_n\stackrel{p}{\rightarrow} 0$, then there exists a
local minimizer $\tilde{\mathbf{V}}_n$ of $Q(\mathbf{V}; \mathbf
{M}_n)$ subject
to $\mathbf{V}^T\mathbf{N}_n\mathbf{V}=\mathbf{I}_d$, so that
\[
D(\tilde{\mathbf{V}}_n,\mathbf{V}_0)=O_p(n^{-1/2}).
\]
\end{theorem}

It is clear from Theorem~\ref{theo1} that by choosing the $\theta_i$'s
properly, there exists a root-$n$ consistent CISE. The next
transition theorem states an oracle-like property of CISE.
\begin{theorem}\label{theo2}
If Assumptions~\ref{Assum1} and~\ref{Assum2} hold, $\lambda
_d>\lambda_{d+1}$,
$\sqrt{n}a_n\stackrel{p}{\rightarrow} 0$ and
$\sqrt{n}b_n\stackrel{p}{\rightarrow} \infty$, then the root-n
consistent local minimizer $\tilde{\mathbf{V}}_n$ in Theorem \ref
{theo1} must
satisfy:
\begin{longlist}
\item $\Pr(\tilde{\mathbf{V}}_{n(p-q)}=\mathbf{0})\rightarrow1$,

\item
$\sqrt{n}D({\tilde{\mathbf{V}}_{n(q)}},\widehat{\mathbf
{V}}_{n(O)})=o_p(1)$, where $\widehat{\mathbf{V}}_{n(O)}$ is the
minimizer of $Q(\mathbf{V}$; $\mathbf{M}_{n(q)})$ subject to $\mathbf
{V}^T\mathbf{N}_{n(q)}\mathbf{V}=\mathbf{I}_d$.
\end{longlist}
\end{theorem}

Theorem~\ref{theo2}(i) states that with probability tending to 1, all
of the
zero row of $\mathbf{V}_0$ must be estimated as $\mathbf{0}$. Theorem
\ref{theo2}(ii)
tells us that there exist a local minimizer $\tilde{\mathbf{V}}_{n}$ so
that the difference between its nonzero submatrix ${\tilde{\mathbf
{V}}_{n(q)}}$ and $\widehat{\mathbf{V}}_{n(O)}$ is of order
$o_p(n^{-1/2})$. That is to say, we have the result that
$\sqrt{n}D({\tilde{\mathbf{V}}_{n(q)}},\mathbf{V}_{0(q)})$ has the same
asymptotic distribution as $\sqrt{n}D(\widehat{\mathbf
{V}}_{n(O)},\mathbf{V}_{0(q)})$. With respect to the asymptotic
distribution of
$\widehat{\mathbf{V}}_{n(O)}$, there seems to be no general result in
the literature because different specifications on $\mathbf{M}_{n(q)}$
and $\mathbf{N}_{n(q)}$ yield different asymptotic distributions. This
is not of great interest here and we refer to \citet{ZN95},
\citet{LZ07} and the references therein.

The second part of Theorem~\ref{theo2} is actually valid in a generalized
sense. The OSDRE in the exact oracle property, denoted as $\dot\mathbf
{V}_{n(O)}$, is obtained by using the $q\times q$ $\mathbf{M}_n$ and
$\mathbf{N}_n$ formed with the first $q$ variables (denoted as
$\mathbf{M}_{n(O)}$ and $\mathbf{N}_{n(O)}$). Usually, $\mathbf
{N}_{n(O)}=\mathbf{N}_{n(q)}$. From the definition, it is
straightforward to see that
$\mathbf{M}_{n(O)}=\mathbf{M}_{n(q)}$ for the PCA, SIR and PFC methods.
Thus, in these cases, Theorem~\ref{theo2} establishes the exact oracle
property. We conjecture that $\mathbf{M}_{n(O)}$ should be very close
to $\mathbf{M}_{n(q)}$ for any SDR method that satisfies Assumptions
\ref{Assum1} and~\ref{Assum2}.
From the proof of Theorem~\ref{theo2}(ii), we can conclude that if
%
%
\begin{equation}\label{oc}
\bigl\Vert\mathbf{M}_{n(O)}-\mathbf{M}_{n(q)}\bigr\Vert_s=O_p(a_n),
\end{equation}
the exact oracle property still holds. The next result establishes
that the condition above holds for DR and SAVE under certain
conditions.
\begin{pro}\label{prop3}
Suppose the linearity and constant variance conditions
[Li and Wang (\citeyear{LiWang2007})] hold and $(na_n)^{-1}=O_p(1)$. Then
condition (\ref{oc})
is satisfied for the DR and SAVE methods.
\end{pro}

By this proposition, Theorem~\ref{theo2} and the discussion above, we know
that from asymptotic viewpoints CISE is effective for all of the
commonly used SDR methods. We summarize this major result in the
following theorem.
\begin{theorem}\label{theo3}
Assume that the conditions in Theorem~\ref{theo2} and Proposition \ref
{prop3} hold. Then
the exact oracle property is achieved for the PCA, SIR, PFC, SAVE
and DR methods. That is, $\tilde{\mathbf{V}}_{n}$ has the selection
consistency and
$\sqrt{n}D({\tilde{\mathbf{V}}_{n(q)}},\mathbf{V}_{0(q)})$
has the same asymptotic distribution as $\sqrt{n}D(\dot\mathbf
{V}_{n(O)},\mathbf{V}_{0(q)})$.
\end{theorem}

In this paper, we make no attempt to further analysis general
conditions for the validity of (\ref{oc}), but we think that such
studies certainly warrant future research.

\subsection{Choice of tuning parameters}\label{sec26}

We recommend
using
%
%
\begin{equation}\label{tuning}
\theta_i = \theta\Vert\widehat{\mathbf{v}}_i\Vert_2^{-r},
\end{equation}
where $\widehat{\mathbf{v}}_i$ is the $i$th row vector of the OSDRE
$\widehat{\mathbf{V}}$ defined in (\ref{Vhat}), and $r>0$ is some
pre-specified parameter. Following the suggestions of \citet{Zou2006},
$r=0.5$ is used in both the simulation study and the illustration in
Section~\ref{sec4}. Such a strategy effectively transforms the original
$p$-dimensional tuning parameter selection problem into a univariate
one. By Lemma~\ref{lema2} in the \hyperref[app]{Appendix}, $\widehat
{\mathbf
{v}}_i$ is root-$n$
consistent. Thus, it is easily to verify that the tuning parameter
defined in (\ref{tuning}) satisfies the conditions on $a_n$ and
$b_n$ needed by Theorem~\ref{theo2} as long as $\sqrt{n}\theta
\rightarrow0$
and ${n}^{(1+r)/2}\theta\rightarrow\infty$. Hence, it suffices to
select $\theta\in[0,+\infty)$ only.

To choose the tuning parameter $\theta$, we use the following
criterion which has a form similar to ones used by \citet{Li2007} and
\citet{LengWang2009}:
\[
-\tr(\tilde{\mathbf{V}}_{\theta}^T\mathbf{M}_n\tilde{\mathbf
{V}}_{\theta})+\gamma\cdot\mathrm{df}_{\theta},
\]
where $\tilde{\mathbf{V}}_{\theta}$ denotes the solution for
$\mathbf{V}$
given $\theta$, $\mathrm{df}_{\theta}$ denotes the effective number of
parameters, and $\gamma=2/n$ for AIC-type and $\gamma=\log(n)/n$ for
BIC-type criteria. Following the discussion of \citet{Li2007}, we
estimate $\mathrm{df}_{\theta}$ by $(p_\theta-d)\cdot d$ where
$p_\theta$ denotes the number of nonzero rows of $\tilde{\mathbf
{V}}_{\theta}$ because we need $(p_\theta-d)\cdot d$ parameters to
describe a $d$-dimensional Grassmann manifold in ${\mathbb
R}^{p_\theta}$ [\citet{EAS98}].

\section{Simulation studies}\label{sec3}

We report the results of four simulation studies in this section,
three of which were
conducted using forward regression models and one was conducted
using an inverse regression model. We compared our method with the
$C^3$ method [\citet{ZhouHe2008}] and the SSIR method
[\citet{NCT05}]. BIC and RIC [\citet{ST02}] were used in SSIR
to select
the tuning parameters, and two $\alpha$ levels (0.01 and 0.005) were
used in the $C^3$ method. We used SIR and PFC to generate $\mathbf
{M}_n$ and $\mathbf{N}_n$ for CISE selection. For these methods, denoted
CIS-SIR and CIS-PFC, we report only the results using the BIC
criterion to select tuning parameters as we tend to believe that BIC
has consistency property. Unreported simulations using the RIC
criterion show slightly better performance in some cases though.

In each study, we generated 2500 datasets with the sample size
$n=60$ and $n=120$. For the $C^3$ method, the quadratic spline with
four internal knots was used, as suggested by \citet{ZhouHe2008}.
Six slices were used for the SSIR method. We calculated $\mathbf{M}_n$
in the PFC model setting using $f(y)=(|y|,y,y^2)^T$ for all
simulation studies.

We used three summary statistics---$r_1$, $r_2$ and $r_3$---to
assess how well the methods select variables: $r_1$ is the average
fraction of nonzero rows of $\tilde{\mathbf{V}}$ associated
with~relevant predictors; $r_2$ is the average fraction of zero rows of
$\tilde{\mathbf{V}}$ associated with irrelevant predictors; and $r_3$ is
the fraction of runs in which the methods select both relevant and
irrelevant predictors exactly right.
%
\begin{Study}\label{Study1}
\[
y = x_1 + x_2 + x_3 + 0.5\epsilon,
\]
where $\epsilon\sim N(0,1)$, ${\mathbf x} = (x_1,\ldots,x_{24})^T
\sim
N(0,\bolds\Sigma)$ with
$\Sigma_{ij} = 0.5^{|i-j|}$ for $1 \leq i,j \leq24$, and ${\mathbf x}$
and $\epsilon$ are independent. In this study, the central subspace
is spanned by the direction $ {\bolds\beta}_1=(1,1,1,0,\ldots,0)^T$
with twenty-one zero coefficients.
\end{Study}
%
%
\begin{Study}\label{Study2}
\[
y = x_1 + x_2 + x_3 + 2\epsilon,
\]
where $\epsilon\sim N(0,1)$, ${\mathbf x} = (x_1,\ldots,x_{24})^T
\sim
N(0,\bolds\Sigma)$ with
$\Sigma_{ij} = 0.5^{|i-j|}$ for $1 \leq i,j \leq24$, and $x$ and
$\epsilon$ are independent. In this study, the central subspace is
spanned by the direction $ {\bolds\beta}_1=(1,1,1,0,\ldots,0)^T$ with
twenty-one zero coefficients. In short, this study was identical to
the first, except the error was increased by a factor of $4$.
\end{Study}
\begin{Study}\label{Study3}
\[
y = x_1/\{0.5 + (x_2 + 1.5)^2\} + 0.2\epsilon,
\]
where $\epsilon\sim N(0,1)$, ${\mathbf x} = (x_1,\ldots,x_{24})^T
\sim
N(0,\bolds\Sigma)$ with
$\Sigma_{ij} = 0.5^{|i-j|}$ for $1 \leq i,j \leq24$, and $x$ and
$\epsilon$ are independent. In this study, the central subspace is
spanned by the directions ${\bolds\beta}_1=(1,0,\ldots,0)^T$ and
$\bolds
\beta_2=(0,1,\ldots,0)^T$.
\end{Study}
\begin{Study}\label{Study4}
\[
{\mathbf x}= \bolds\Gamma(y,y^2)^T + \bolds\Delta^{1/2} \bolds
\epsilon,
\]
where $\bolds\epsilon\sim N(0,\mathbf{I}_{24})$, $y \sim N(0,1)$,
$\Delta
_{ij} = 0.5^{|i-j|}$ for $1 \leq i,j \leq
24$, and $y$ and $\bolds\epsilon$ are independent. The first column of
$\bolds\Gamma$ is $(0.5,0.5,0.5,0.5,0,\ldots,0)^T$ and the second
column of $\bolds\Gamma$ is $(0.5,-0.5,0.5,-0.5,0,\ldots,0)^T$. In this
study, the central subspace is the column space of
$\bolds\Delta^{-1}\bolds\Gamma$.
\end{Study}

The simulation results from these four studies are summarized in Tables~\mbox{\ref{table2}--\ref{table5}}, respectively.
The standard errors of the $r_k$'s, $\sqrt{r_k(1-r_k)}/50$, are
typically less
than 0.01 throughout this section. In Study~\ref{Study1}, the signal-to-noise
ratio is close to 5 (the ratio of the stand deviation of $ x_1 + x_2
+ x_3$ to 0.5). Because of the large signal-to-noise ratio, all the
considered methods show very good
performance, but CIS-SIR, CIS-PFC and $C^3$ perform slightly
better than SSIR. In Study~\ref{Study2}, we decreased the
signal-to-noise ratio to about 1.2 and now
CIS-SIR and CIS-PFC perform much better than $C^3$ and SSIR.
In both Studies~\ref{Study3}
and~\ref{Study4}, CISE is generally superior to the other two methods,
especially for CIS-PFC and the rate $r_3$. It should be pointed
out that the superiority of CISE becomes more significant when $n$
gets larger. When $n=120$, $C^3$ still cannot perform exact
identifications well, while SSIR rarely identifies all
relevant and irrelevant variables correctly.\looseness=1

%
%
\begin{table}
\caption{Summary of Study \protect\ref{Study1}}\label{table2}
\begin{tabular*}{\tablewidth}{@{\extracolsep{\fill}}lcccccc@{}}
\hline
\multicolumn{1}{@{}r}{\textbf{Method:}}
& \multicolumn{1}{c}{\textbf{CIS-SIR}} & \multicolumn{1}{c}{\textbf{CIS-PFC}} &
\multicolumn{2}{c}{$\bolds{C^3}$} & \multicolumn{2}{c@{}}{\textbf{SSIR}}\\[-4pt]
& \multicolumn{1}{c}{\hrulefill} & \multicolumn{1}{c}{\hrulefill}
& \multicolumn{2}{c}{\hrulefill} & \multicolumn{2}{c@{}}{\hrulefill}
\\
\multicolumn{1}{@{}r}{\textbf{Criterion:}} &
\multicolumn{1}{c}{\textbf{BIC}} & \multicolumn{1}{c}{\textbf{BIC}} & \multicolumn{1}{c}{$\bolds{\alpha=0.01}$}
& \multicolumn{1}{c}{$\bolds{\alpha=0.005}$} & \multicolumn{1}{c}{\textbf{BIC}} & \multicolumn{1}{c@{}}{\textbf{RIC}} \\
\hline
Sample size & \multicolumn{6}{c@{}}{$n=60$} \\[3pt]
\quad$r_1$ & 0.991 & 1.000 & 1.000 & 1.000 & 0.993 & 0.974 \\
\quad$r_2$ & 0.999 & 1.000 & 0.999 & 0.999 & 0.997 & 0.999 \\
\quad$r_3$ & 0.970 & 1.000 & 0.978 & 0.991 & 0.939 & 0.914\\[3pt]
Sample size & \multicolumn{6}{c@{}}{$n=120$} \\[3pt]
\quad$r_1$ & 1.000 & 1.000 & 1.000 & 1.000 & 1.000 & 1.000 \\
\quad$r_2$ & 1.000 & 1.000 & 1.000 & 1.000 & 0.999 & 1.000 \\
\quad$r_3$ & 1.000 & 1.000 & 1.000 & 1.000 & 0.994 & 1.000
\\
\hline
\end{tabular*}\vspace*{-1.5pt}
\end{table}

%
%
\begin{table}
\caption{Summary of Study \protect\ref{Study2}}
\label{table3}
\begin{tabular*}{\tablewidth}{@{\extracolsep{\fill}}lcccccc@{}}
\hline
\multicolumn{1}{@{}r}{\textbf{Method:}}
& \multicolumn{1}{c}{\textbf{CIS-SIR}} & \multicolumn{1}{c}{\textbf{CIS-PFC}} &
\multicolumn{2}{c}{$\bolds{C^3}$} & \multicolumn{2}{c@{}}{\textbf{SSIR}}\\[-4pt]
& \multicolumn{1}{c}{\hrulefill} & \multicolumn{1}{c}{\hrulefill}
& \multicolumn{2}{c}{\hrulefill} & \multicolumn{2}{c@{}}{\hrulefill}
\\
\multicolumn{1}{@{}r}{\textbf{Criterion:}} &
\multicolumn{1}{c}{\textbf{BIC}} & \multicolumn{1}{c}{\textbf{BIC}} & \multicolumn{1}{c}{$\bolds{\alpha=0.01}$}
& \multicolumn{1}{c}{$\bolds{\alpha=0.005}$} & \multicolumn{1}{c}{\textbf{BIC}} & \multicolumn{1}{c@{}}{\textbf{RIC}} \\
\hline
Sample size & \multicolumn{6}{c@{}}{$n=60$} \\[3pt]
\quad$r_1$ & 0.713 & 0.795 & 0.583 & 0.565 & 0.770 & 0.706 \\
\quad$r_2$ & 0.988 & 0.992 & 0.998 & 0.998 & 0.881 & 0.939 \\
\quad$r_3$ & 0.233 & 0.399 & 0.075 & 0.080 & 0.058 & 0.104 \\
[3pt]
Sample size & \multicolumn{6}{c@{}}{$n=120$} \\[3pt]
\quad$r_1$ & 0.909 & 0.951 & 0.669 & 0.615 & 0.973 & 0.930 \\
\quad$r_2$ & 0.998 & 0.998 & 1.000 & 1.000 & 0.928 & 0.981 \\
\quad$r_3$ & 0.694 & 0.827 & 0.209 & 0.131 & 0.244 & 0.554
\\\hline
\end{tabular*}\vspace*{-1.5pt}
\end{table}

%
%
\begin{table}
\caption{Summary of Study \protect\ref{Study3}}\label{table4}
\begin{tabular*}{\tablewidth}{@{\extracolsep{\fill}}lcccccc@{}}
\hline
\multicolumn{1}{@{}r}{\textbf{Method:}}
& \multicolumn{1}{c}{\textbf{CIS-SIR}} & \multicolumn{1}{c}{\textbf{CIS-PFC}} &
\multicolumn{2}{c}{$\bolds{C^3}$} & \multicolumn{2}{c@{}}{\textbf{SSIR}}\\[-4pt]
& \multicolumn{1}{c}{\hrulefill} & \multicolumn{1}{c}{\hrulefill}
& \multicolumn{2}{c}{\hrulefill} & \multicolumn{2}{c@{}}{\hrulefill}
\\
\multicolumn{1}{@{}r}{\textbf{Criterion:}} &
\multicolumn{1}{c}{\textbf{BIC}} & \multicolumn{1}{c}{\textbf{BIC}} & \multicolumn{1}{c}{$\bolds{\alpha=0.01}$}
& \multicolumn{1}{c}{$\bolds{\alpha=0.005}$} & \multicolumn{1}{c}{\textbf{BIC}} & \multicolumn{1}{c@{}}{\textbf{RIC}} \\
\hline
Sample size & \multicolumn{6}{c@{}}{$n=60$} \\[3pt]
\quad$r_1$ & 0.789 & 0.906 & 0.770 & 0.742 & 0.934 & 0.888 \\
\quad$r_2$ & 0.965 & 0.979 & 0.948 & 0.955 & 0.633 & 0.828 \\
\quad$r_3$ & 0.344 & 0.588 & 0.229 & 0.226 & 0.000 & 0.004 \\
[3pt]
Sample size & \multicolumn{6}{c@{}}{$n=120$} \\[3pt]
\quad$r_1$ & 0.948 & 0.995 & 0.839 & 0.781 & 0.994 & 0.983 \\
\quad$r_2$ & 0.992 & 0.998 & 0.956 & 0.963 & 0.664 & 0.865 \\
\quad$r_3$ & 0.838 & 0.973 & 0.309 & 0.245 & 0.001 & 0.027
\\\hline
\end{tabular*}\vspace*{-1.5pt}
\end{table}

%
%
\begin{table}
\caption{Summary of Study \protect\ref{Study4}}\label{table5}
\begin{tabular*}{\tablewidth}{@{\extracolsep{\fill}}lcccccc@{}}
\hline
\multicolumn{1}{@{}r}{\textbf{Method:}}
& \multicolumn{1}{c}{\textbf{CIS-SIR}} & \multicolumn{1}{c}{\textbf{CIS-PFC}} &
\multicolumn{2}{c}{$\bolds{C^3}$} & \multicolumn{2}{c@{}}{\textbf{SSIR}}\\[-4pt]
& \multicolumn{1}{c}{\hrulefill} & \multicolumn{1}{c}{\hrulefill}
& \multicolumn{2}{c}{\hrulefill} & \multicolumn{2}{c@{}}{\hrulefill}
\\
\multicolumn{1}{@{}r}{\textbf{Criterion:}} &
\multicolumn{1}{c}{\textbf{BIC}} & \multicolumn{1}{c}{\textbf{BIC}} & \multicolumn{1}{c}{$\bolds{\alpha=0.01}$}
& \multicolumn{1}{c}{$\bolds{\alpha=0.005}$} & \multicolumn{1}{c}{\textbf{BIC}} & \multicolumn{1}{c@{}}{\textbf{RIC}} \\
\hline
Sample size & \multicolumn{6}{c@{}}{$n=60$} \\[3pt]
\quad$r_1$ & 0.676 & 0.817 & 0.670 & 0.643 & 0.871 & 0.776 \\
\quad$r_2$ & 0.968 & 0.989 & 0.956 & 0.958 & 0.641 & 0.832 \\
\quad$r_3$ & 0.069 & 0.327 & 0.022 & 0.029 & 0.000 & 0.000 \\
[3pt]
Sample size & \multicolumn{6}{c@{}}{$n=120$} \\[3pt]
\quad$r_1$ & 0.805 & 0.928 & 0.828 & 0.809 & 0.988 & 0.964 \\
\quad$r_2$ & 0.993 & 0.998 & 0.967 & 0.969 & 0.696 & 0.890 \\
\quad$r_3$ & 0.299 & 0.687 & 0.147 & 0.178 & 0.000 & 0.000
\\\hline
\end{tabular*}
\end{table}

While both CISE and $C^3$ have the oracle property,
they differ in many aspects. CISE is a unified
method that can be applied to many popular sufficient dimension
reduction methods, including PCA, PFC, SIR, SAVE and DR. On the other hand,
$C^3$ is based on one specified sufficient dimension reduction
method, canonical correlation [\citet{Fungetal02}]. We regard $r_{3}$,
the estimated probability all relevant and irrelevant variables are
identified correctly, as the most important aspect of a method. On
that measure CISE typically dominates $C^{3}$. There was only one
case (Table~\ref{table1}, $n=60$) in which $C^{3}$ did slightly better than
CISE.
Additionally, CISE seems conceptually simpler and is
easily implemented.

\section{Boston housing data}\label{sec4}
\subsection{Variable screening}
We applied our method to the Boston housing data, which has been
widely studied in the literature. The Boston housing data contains
506 observations, and can be downloaded from the web site
\url{http://lib.stat.cmu.edu/datasets/boston\_corrected.txt}. The
re-\break sponse~variable $y$ is the median value of owner-occupied homes in
each of the 506 census tracts in the Boston Standard Metropolitan
Statistical Areas. The 13 predictor variables are per capita crime
rate by town ($x_1$); proportion of residential land zoned for lots
over 25,000 sq.ft ($x_2$); proportion of nonretail business acres
per town ($x_3$); Charles River dummy variable ($x_4$); nitric
oxides concentration ($x_5$); average number of rooms per dwelling
($x_6$); proportion of owner-occupied units built prior to 1940
($x_7$); weighted distances to five Boston employment centers
($x_8$); index of accessibility to radial highways ($x_9$);
full-value property-tax rate ($x_{10}$); pupil--teacher ratio by town
($x_{11}$); proportion of blacks by town ($x_{12}$); percentage of
lower status of the population~($x_{13}$).

%
%
\begin{table}
\tabcolsep=0pt
\caption{Estimated bases of the central subspace in Boston housing data}
\label{table6}
{\fontsize{8.66pt}{10.36pt}\selectfont{
\begin{tabular*}{\tablewidth}{@{\extracolsep{\fill}}lk{6.7}k{6.7}k{6.7}k{6.7}k{6.7}@{\hspace*{-2pt}}}
\hline
\textbf{Method:} & \multicolumn{1}{c}{\textbf{CIS-SIR}} &
\multicolumn{1}{c}{\textbf{CIS-PFC}} & \multicolumn{1}{c}{$\bolds{C^3}$} &
\multicolumn{1}{c}{\textbf{SSIR-BIC}}& \multicolumn{1}{c@{}}{\textbf{SSIR-RIC}}\\
\hline
$x_1$  & 0, \mbox{ }0 & 0, \mbox{ }0 & 0, \mbox{ }0 & \mbox{$-$}0.050, \mbox{ }\mbox{$-$}0.131 & \mbox{$-$}0.041, \mbox{ }\mbox{$-$}0.123\\
$x_2$ & \mbox{$-$}0.004, \mbox{ }\mbox{$-$}0.047 & 0, \mbox{ }0 & 0, \mbox{ }0 &\mbox{$-$}0.001,\mbox{ }\mbox{$-$}0.002 & \mbox{$-$}0.001, \mbox{ }\mbox{$-$}0.001\\
$x_3$ & 0, \mbox{ }0 &0, \mbox{ }0 & 0, \mbox{ }0 &  0.001, \mbox{ } 0.005& 0, \mbox{ }0\\
$x_4$ & 0, \mbox{ }0 &0, \mbox{ }0 & 0, \mbox{ }0 & \mbox{$-$}0.033, \mbox{ }0.020 & 0, \mbox{ }0 \\
$x_5$ & 0, \mbox{ }0 &0, \mbox{ }0 & 0, \mbox{ }0&  0.719, \mbox{ }\mbox{$-$}0.882 &  0.543, \mbox{ }\mbox{$-$}0.765 \\
$x_6$ & \mbox{$-$}0.999, \mbox{ }0.034 &\mbox{$-$}0.999, \mbox{ }0.034 & 0.962, \mbox{ }\mbox{$-$}0.645 &  \mbox{$-$}0.684, \mbox{ }\mbox{$-$}0.448 & \mbox{$-$}0.834, \mbox{ }\mbox{$-$}0.627\\
$x_7$ & \mbox{$-$}0.008, \mbox{ }\mbox{$-$}0.139 &\mbox{$-$}0.003, \mbox{ }\mbox{$-$}0.077 & \mbox{$-$}0.174, \mbox{ }\mbox{$-$}0.096&  0.006, \mbox{ }\mbox{$-$}0.001& 0.005, \mbox{ }\mbox{$-$}0.001\\
$x_8$ & 0, \mbox{ }0 &0, \mbox{ }0 & 0, \mbox{ }0 & 0.082, \mbox{ }\mbox{$-$}0.012 & 0.060, \mbox{ }\mbox{$-$}0.010\\
$x_9$ & 0, \mbox{ }0 &0, \mbox{ }0 & 0, \mbox{ }0 &  \mbox{$-$}0.019, \mbox{ }0.035 & \mbox{$-$}0.016, \mbox{ }0.033 \\
$x_{10}$ & \mbox{$-$}0.001, \mbox{ }\mbox{$-$}0.01 & \mbox{$-$}0.002, \mbox{ }\mbox{$-$}0.035 & \mbox{$-$}0.166, \mbox{ }0 & 0.001, \mbox{ }\mbox{$-$}0.001& 0.001, \mbox{ }\mbox{$-$}0.001\\
$x_{11}$ & 0.021, \mbox{ }\mbox{$-$}0.361 &0.018, \mbox{ }\mbox{$-$}0.280 & \mbox{$-$}0.126, \mbox{ }0 &  0.058, \mbox{ }\mbox{$-$}0.033& 0.055, \mbox{ }\mbox{$-$}0.036\\
$x_{12}$ & 0.001, \mbox{ }0.011 &0.002, \mbox{ }0.035 & 0, \mbox{ }0 &  \mbox{$-$}0.000, \mbox{ }0.000 &  0, \mbox{ }0 \\
$x_{13}$ & \mbox{$-$}0.044, \mbox{ }\mbox{$-$}0.920 &\mbox{$-$}0.040, \mbox{ }\mbox{$-$}0.955 & 0, \mbox{ }\mbox{$-$}0.758 & 0.014, \mbox{ }\mbox{$-$}0.043& 0.017, \mbox{ }\mbox{$-$}0.059\\
\hline
\end{tabular*}
}}
\end{table}

%
%
\begin{table}[b]
\caption{Variable selection in bootstrapping Boston housing data}
\label{table7}
\begin{tabular*}{\tablewidth}{@{\extracolsep{\fill}}lcccc@{}}
\hline
\textbf{Method:} & \multicolumn{1}{c}{\textbf{CIS-SIR}} &
\multicolumn{1}{c}{\textbf{CIS-PFC}} & \multicolumn{1}{c}{\textbf{SSIR-BIC}}&
\multicolumn{1}{c@{}}{\textbf{SSIR-RIC}}\\
\hline
$r_1$ & 0.947 & 0.962 & 0.963 & 0.877 \\
$r_2$ & 0.969 & 0.980 & 0.780 & 0.952 \\
$r_3$ & 0.550 & 0.672 & 0.118 & 0.264 \\
\hline
\end{tabular*}
\end{table}

Previous studies suggested that we remove those observation with crime
rate greater than 3.2, as a few predictors remain constant except for
3 observations in this case [\citet{L91}]. So we used the 374
observations with crime rate smaller than 3.2 in this analysis. All
the methods considered in Section~\ref{sec3} were applied to this dataset.
Scatter-plotting of each predictor against $y$, we concluded that it
would be sufficient to use $\mathbf{f} = (\sqrt{y},y,y^2)^{T}$ in the
PFC model. Since PFC is a scale-invariant method, we did not
standardize the data as many other methods do. Similar to the previous
studies in the literature, we pick up two directions to estimate the
central subspace. The estimated bases of the central subspace for
all the considered methods are summarized in Table~\ref{table6}.

The coefficients in Table~\ref{table6} from CIS-SIR, CIS-PFC and SSIR
are based on the original dataset, while the coefficients of $C^3$
is based on a data-specific weighted version [\citet{ZhouHe2008}].
As suggested by CIS-PFC, explanatory variables $x_6$, $x_7$, $x_{10}$,
$x_{11}$, $x_{12}$ and $x_{13}$ would be important in explaining
$y$.

\subsection{Bootstrap study}

In Table~\ref{table7}, we used the bootstrap to assess the accuracy of variable
selection for all methods except $C^3$, as it is not clear how the
weighting procedure used by Zhou and He should be automated. Without
weighting we encountered serious convergence problems in the $C^{3}$
algorithm. This bootstrap study can be considered as another
simulation study.

The bootstrap procedure was conducted as follows. First, we
randomly chose with replacement 374 observations for $y$ jointly
with $x_6$, $x_7$, $x_{10}$, $x_{11}$, $x_{12}$ and $x_{13}$.
Secondly, we separately randomly selected 374 observations for
$x_1$, $x_2$, $x_3$, $x_4$, $x_5$, $x_8$ and $x_9$. Then we combine
them to make one complete bootstrap dataset. In this way, we mimic
the results of the analysis of original data, forcing $x_1$, $x_2$,
$x_3$, $x_4$, $x_5$, $x_8$ and $x_9$ to be irrelevant. This
procedure was repeated 2500 times. The resulting rates $r_1$, $r_2$
and $r_3$ are shown in Table~\ref{table7}. 
The results show a pattern similar to those in simulation
studies and again CISE performed quite well.

\section{Discussion}\label{sec5}

The establishment of the oracle property in this paper~takes
advantage of the simple trace form of the objective function:
$-\tr(\mathbf{V}^T\mathbf{M}_n\mathbf{V})$. However we believe that
the proof
in the \hyperref[app]{Appendix} can be extended to more general
objective functions.
Moreover, it is also of great interests to see whether CISE and its
oracle property are still valid in high-dimensional settings in
which $p>n$.

We have seen that $\mathbf{N}_n$ usually takes the form of the marginal
sample covariance matrix of ${\mathbf x}$, while $\mathbf{M}_n$ depends
on the specific method. In practice, how to choose $\mathbf{M}_n$ for
variable selection is an
important issue and merits thorough investigation. In addition, it
is well demonstrated that for the multiple regression model, the BIC
criterion tends to identify the true sparse model well if the true
model is included in the candidate set [\citet{WLT07}]. The
consistency of the BIC criterion proposed in Section~\ref{sec26} deserves
further study as well.

\begin{appendix}\label{app}
\section*{Appendix}

Throughout this section, we will use the following notation for ease
of exposition. $Q(\bolds\Gamma;\mathbf{G}_n,\mathbf{N}_n):=-\tr
(\bolds\Gamma^T\mathbf{G}_n\bolds\Gamma)+\rho(\mathbf
{N}_n^{-1/2}\bolds\Gamma)$ denotes the constrained objective function in
the minimization problem (\ref{GammaTilde}). Unless otherwise
stated, we also use the generic notation $Q(\bolds\Gamma)$ or
$Q(\mathbf{V})$ to represent the function $Q(\bolds\Gamma;\mathbf
{G}_n,\mathbf{N}_n)$ or
$Q(\mathbf{V}; \mathbf{M}_n)$ for abbreviation, which should not
cause any
confusion. $\mathbf{1}_i$ denotes a row vector with one in the $i$th
position and zero in the others.
\begin{pf*}{Proof of Proposition~\ref{prop1}}
\citet{Cook07} has shown that the
maximum likelihood estimator of
$\operatorname{span}(\bolds\Delta^{-1}\bolds\Gamma)$ in the general PFC model
equals the span of $\{\mathbf{e}_1,\ldots,\mathbf{e}_d\}$, where
$\mathbf{e}_i={\bolds\Sigma}_n^{-1/2}\mathbf{r}_i$ and $\mathbf
{r}_i$ is the $i$th
eigenvector of
${\bolds\Sigma}_n^{-1/2}\widehat{\bolds\Sigma}_{\mathrm
{fit}}{\bolds\Sigma}_n^{-1/2}$
corresponding to the eigenvalue $k_i$. Consequently, we have
\[
\widehat{\bolds\Sigma}_{\mathrm{fit}} \mathbf{e}_i = k_i {\bolds
\Sigma}_n
\mathbf{e}_i.
\]
It follows that $\mathbf{M}_n=\widehat{\bolds\Sigma
}_{\mathrm{fit}}$ and $\mathbf{N}_n
={\bolds\Sigma}_n$.
\end{pf*}

In order to prove the theorems, we first state a
few necessary lemmas. For notation convenience, we need the
following additional definitions. Define the Stiefel manifold
$\operatorname{St}(p,d)$ as
\[
\operatorname{St}(p,d)=\{\bolds\Gamma\in{\mathbb R}^{p\times
d}\dvtx\bolds\Gamma^T\bolds\Gamma=\mathbf{I}_d\}.
\]
Denotes $\lfloor\bolds\Gamma\rfloor$ as the subspace spanned by the
columns of $\bolds\Gamma$, then $\lfloor\bolds\Gamma\rfloor\in \operatorname{Gr}(p,d)$
where $\operatorname{Gr}(p,d)$ stands for the Grassmann manifold.\vadjust{\goodbreak} The projection
operator $R\dvtx{\mathbb R}^{p\times d}\rightarrow \operatorname{St}(p,d)$ onto the
Stiefel manifold $\operatorname{St}(p,d)$ is defined to~be
\[
R(\bolds\Gamma)={\mathop{\arg\min}_{\mathbf{W}\in \operatorname{St}(p,d)}}\Vert
\bolds\Gamma-\mathbf{W}\Vert^2_s.
\]
The tangent space $T_{\bolds\Gamma}(p,d)$ of $\bolds\Gamma\in
\operatorname{St}(p,d)$ is
defined by
%
%
\begin{eqnarray}\label{eqA1}
T_{\bolds\Gamma}(p,d)&=&\bigl\{\mathbf{Z}\in{\mathbb R}^{p\times d}\dvtx
\mathbf{Z}=\bolds\Gamma\mathbf{A}+\bolds\Gamma_{\bot}\mathbf
{B},\nonumber\\[-8pt]\\[-8pt]
&&\hspace*{4.5pt} \mathbf{A}\in{\mathbb
R}^{d\times d}, \mathbf{A}+\mathbf{A}^T=0, \mathbf{B}\in{\mathbb
R}^{(p-d)\times d}\bigr\},\nonumber
\end{eqnarray}
where $\bolds\Gamma_{\bot}\in{\mathbb R}^{p\times{(p-d)}}$ is the
complement of $\bolds\Gamma$ satisfies $[\bolds\Gamma
\bolds\Gamma_{\bot}]^T[\bolds\Gamma\bolds\Gamma_{\bot}]=\mathbf{I}_p$.
%
%
\begin{lemma}\label{lema1}
If $\mathbf{Z}\in T_{\bolds\Gamma}(p,d), {\bolds
\Gamma}\in
\operatorname{St}(p,d)$, we have:
\begin{longlist}
\item For any symmetric matrix $\mathbf{C}\in{\mathbb
R}^{d\times
d}$, $\tr(\mathbf{Z}^T\bolds\Gamma\mathbf{C})=0$.

\item $R(\bolds\Gamma+t\mathbf{Z})=\bolds\Gamma+t\mathbf
{Z}-(1/2)t^2\bolds\Gamma
\mathbf{Z}^T\mathbf{Z}+O(t^3)$.
\end{longlist}
\end{lemma}

This lemma comes from Lemma 10 and Proposition 12 of \citet{Manton2002}.

\begin{lemma}\label{lema2}
Under conditions in Theorem~\ref{theo1}, we have
\[
D({\widehat{\bolds\Gamma}},{\bolds\Gamma_0})=O_p(n^{-1/2}),
\]
where ${\bolds\Gamma_0}$ denotes any minimizer of (\ref{GammaHat}) when
$\mathbf{G}_n$ is taken as the population matrix~$\mathbf{G}$.
\end{lemma}

This lemma can be proved in a similar fashion to the proof of
Theorem 1 and hence omitted here.
\begin{pf*}{Proof of Theorem~\ref{theo1}}
Clearly, to prove this theorem is
equivalent to show there exists a local minimizer $\tilde{\bolds
\Gamma}_n$ of $Q({\bolds\Gamma}; \mathbf{G}_n,\mathbf{N}_n)$
subject to
${\bolds\Gamma}^T{\bolds\Gamma}=\mathbf{I}_d$, so that
\[
D(\tilde{\bolds\Gamma}_n,\bolds{\Gamma}_{0})=O_p(n^{-1/2}).
\]
Denote $\bolds\Gamma_{*}$ as an orthonormal basis matrix of the
subspace spanned by the columns of $\mathbf{N}_n^{1/2}\mathbf{V}_0$. Thus,
there exists a positive-definite matrix $\mathbf{O}\in{\mathbb
R}^{d\times d}$ so that $\bolds\Gamma_{*}=\mathbf{N}_n^{1/2}\mathbf
{V}_0\mathbf{O}$. By Assumption~\ref{Assum2} and $\mathbf
{V}_0^{T}\mathbf{N}\mathbf{V}_0=\mathbf{I}_d$,
we have
\[
\mathbf{O}^{T}\mathbf{O}=\mathbf{I}_d+O_p(n^{-1/2}).
\]
Note that
$\bolds\Gamma_{0} = \mathbf{N}^{1/2}\mathbf{V}_0$, and thus it is equivalent
to show that
\[
D(\tilde{\bolds\Gamma}_n,\bolds{\Gamma}_{*})=O_p(n^{-1/2}),
\]
since $D(\bolds{\Gamma}_{*},{\bolds\Gamma}_0)=O_p(n^{-1/2})$ and
$D(\cdot,\cdot)$ satisfies the triangle inequality.

To ease demonstration, we need define the concept of the
neighborhood of $\lfloor\bolds\Gamma_{*}\rfloor$. For an arbitrary
matrix $\mathbf{W}\in{\mathbb R}^{p\times d}$ and scaler $\delta\in
{\mathbb R}$, the perturbed point around ${\bolds\Gamma_{*}}$ in
Stiefel manifold can be expressed by $R({\bolds\Gamma_{*}}+\delta
\mathbf{W})$. The perturbed point around $\lfloor\bolds\Gamma
_{*}\rfloor$ in
Grassmann manifold can be expressed by $\lfloor
R({\bolds\Gamma_{*}}+\delta\mathbf{W})\rfloor$. According to Lemma
8 of
\citet{Manton2002}, $\mathbf{W}$ can be uniquely decomposed as
\[
\mathbf{W} =
\bolds\Gamma_{*} \mathbf{A}+{\bolds\Gamma_{*}}_{\bot}\mathbf{B} +
\bolds\Gamma_{*}\mathbf{C},
\]
where $\mathbf{A} \in{\mathbb R}^{d\times d}$
is a skew-symmetric matrix, $\mathbf{B} \in{\mathbb R}^{(p-d)\times
d}$ is an arbitrary matrix, and $\mathbf{C} \in{\mathbb R}^{d\times
d}$ is a symmetric matrix. Let $\mathbf{Z} =\bolds\Gamma_{*} \mathbf
{A}+{\bolds\Gamma_{*}}_{\bot}\mathbf{B}$. Obviously, $\mathbf
{Z}\in
T_{\bolds\Gamma_{*}}(p,d)$. Henceforth, $\mathbf{Z}$ refers to the
projection of an arbitrary matrix $\mathbf{W}\in{\mathbb R}^{p\times
d}$ onto the tangent space $T_{\bolds\Gamma_{*}}(p,d)$, unless
otherwise stated.

From Proposition 20 of \citet{Manton2002}, it is straightforward
to see
\begin{eqnarray*}
\lfloor R({\bolds\Gamma_{*}}+\delta\mathbf{W})\rfloor&=& \bigl\lfloor
R\bigl({\bolds\Gamma_{*}}+\delta(\bolds\Gamma_{*} \mathbf{A}+{\bolds
\Gamma_{*}}_{\bot}\mathbf{B} + \bolds\Gamma_{*}\mathbf
{C})\bigr)\bigr\rfloor\\
&=& \bigl\lfloor R\bigl({\bolds\Gamma_{*}}\bigl(\mathbf{I}_d + \delta(\mathbf{A}
+\mathbf{C})\bigr)+\delta{\bolds\Gamma_{*}}_{\bot}\mathbf{B}\bigr) \bigr\rfloor
\\
&=& \bigl\lfloor{\bolds\Gamma_{*}}\bigl(\mathbf{I}_d + \delta(\mathbf{A}
+\mathbf{C})\bigr)+\delta{\bolds\Gamma_{*}}_{\bot}\mathbf{B} \bigr\rfloor\\
&=& \bigl\lfloor{\bolds\Gamma_{*}} + \delta{\bolds\Gamma_{*}}_{\bot
}\mathbf{B}\bigl(\mathbf{I}_d + \delta(\mathbf{A} +\mathbf{C})\bigr)^{-1}
\bigr\rfloor\\
&=& \lfloor
R({\bolds\Gamma_{*}}+\delta{\bolds\Gamma_{*}}_{\bot}{\mathbf
B}^{\prime})\rfloor,
\end{eqnarray*}
provided that $\delta$ is sufficiently small so that $\mathbf{I}_d +
\delta(\mathbf{A} +\mathbf{C})$ is a full rank matrix, where
${\mathbf
B}^{\prime} =\mathbf{B} (\mathbf{I}_d + \delta(\mathbf{A} +\mathbf
{C}))^{-1}$.
Since $\mathbf{B} \in{\mathbb R}^{(p-d)\times d}$ is an arbitrary
matrix and we do not need the specific form of $\mathbf{B}$ and
${\mathbf
B}^{\prime}$ in our proof, we only use $\mathbf{B}$ for notation
convenience. This tells us that the movement from
$\lfloor\bolds\Gamma_{*}\rfloor$ in the near neighborhood only depends
on the ${\bolds\Gamma_{*}}_{\bot}\mathbf{B}$. In other words, it suffices
to only consider perturbed points like $R({\bolds\Gamma_{*}}+\delta
\mathbf{Z})$ in the following proofs, where $\Vert\mathbf{B}\Vert
_s=C$ for some
given~$C$. It is worth noting that though our problems essentially
are Grassmann manifold optimization, we prove the theorem in a more
general way, say in Stiefel manifold [using $\mathbf{Z}\in
T_{\bolds\Gamma_{*}}(p,d)$] since the latter has simpler matrix
expressions and thus is more notationally convenient.

For any small $\epsilon$, if we can show that there exits a
sufficiently large constant~$C$, such that
%
%
\begin{equation}\label{eqA2}\qquad
\lim_n \pr\Bigl(\inf_{\mathbf{Z}\in T_{\bolds\Gamma_{*}}(p,d)\dvtx
\Vert\mathbf{B}\Vert_s=C}Q\bigl(R({\bolds\Gamma_{*}}+n^{-
1/2}\mathbf{Z})\bigr)>Q({\bolds\Gamma}_{*})\Bigr)>1-\varepsilon,
\end{equation}
then we can conclude that there exists a local
minimizer $\tilde{\bolds\Gamma}_n$ of $Q({\bolds\Gamma})$ with arbitrarily
large probabilities such that
$\Vert\tilde{\bolds\Gamma}_n-\bolds{\Gamma}_*\Vert
_s=O_p(n^{-1/2})$. This
certainly implies that
$D(\tilde{\bolds\Gamma}_n,\bolds{\Gamma}_*)=O_p(n^{-1/2})$ by Definition
\ref{defi1}.

By using Lemma~\ref{lema1}, for $\mathbf{Z}\in T_{\bolds\Gamma
_{*}}(p,d)$ we have
\begin{eqnarray*}
&&n\bigl\{Q\bigl(R({\bolds\Gamma}_*+n^{-1/2}\mathbf{Z})\bigr)-Q({\bolds\Gamma
}_*)\bigr\}
\\
&&\qquad=\bigl[-\tr(\mathbf{Z}^T\mathbf{G}_n\mathbf{Z})-2\sqrt{n}\tr(\mathbf
{Z}^T{\mathbf{G}}_n{\bolds\Gamma}_*)+\tr(\mathbf{Z}^T\mathbf
{Z}{\bolds\Gamma}_*^T{\mathbf{G}}_n{{\bolds
\Gamma}_*})\bigr]\bigl(1+o_p(1)\bigr)\\
&&\qquad\quad{}+n\sum_{j=1}^p\biggl[\theta_j\biggl\Vert{\mathbf1}_j\mathbf{N}_n^{-1/2}\biggl({\bolds\Gamma}_*+n^{-1/2}\mathbf{Z}-\frac12n^{-1}{\bolds
\Gamma}_*\mathbf{Z}^T\mathbf{Z}\biggr)\biggr\Vert_2\\
&&\qquad\quad\hspace*{158.6pt}{}-\theta_j\Vert{\mathbf
1}_j\mathbf{N}_n^{-1/2}{\bolds\Gamma}_*\Vert_2\biggr]
\bigl(1+o_p(1)\bigr)\\
&&\qquad\geq\bigl[-\tr(\mathbf{Z}^T\mathbf{G}_n\mathbf{Z})-2\sqrt{n}\tr
(\mathbf{Z}^T{\mathbf{G}}_n{\bolds\Gamma}_*)+\tr(\mathbf
{Z}^T\mathbf{Z}{\bolds\Gamma}_*^T{\mathbf{G}}_n{{\bolds
\Gamma}_*})\bigr]\bigl(1+o_p(1)\bigr)\\
&&\qquad\quad{}+n\sum_{j=1}^q\biggl[\theta_j\biggl(\biggl\Vert{\mathbf1}_j\mathbf{N}_n^{-1/2}\biggl({\bolds\Gamma}_*+n^{-1/2}\mathbf{Z}-\frac12n^{-1}{\bolds
\Gamma}_*\mathbf{Z}^T\mathbf{Z}\biggr)\biggr\Vert_2\\
&&\qquad\quad\hspace*{175.5pt}{}-\Vert{\mathbf
1}_j\mathbf{N}_n^{-1/2}{\bolds\Gamma}_*\Vert_2\biggr)\biggr]\bigl(1+o_p(1)\bigr)\\
&&\qquad\geq\bigl[-\tr(\mathbf{Z}^T\mathbf{G}_n\mathbf{Z})+\tr(\mathbf
{Z}^T\mathbf{Z}{\bolds\Gamma}_*^T{\mathbf{G}}_n{{\bolds\Gamma
}_*})-2\sqrt{n}\tr(\mathbf{Z}^T{\mathbf{G}}_n{\bolds\Gamma
}_*)\bigr]\bigl(1+o_p(1)\bigr)\\
&&\qquad\quad{}-\frac12q\bigl(\sqrt{n}a_n\bigr)\\
&&\qquad\quad\hspace*{10.8pt}{}\times\max_j\bigl\{\Vert{\mathbf1}_j\mathbf
{N}_n^{-1/2}\bolds\Gamma_*\Vert_2^{-1}\cdot\bigl\Vert{\mathbf
1}_j\mathbf{N}_n^{-1/2}\bigl(\mathbf{Z}-(1/2)n^{-1/2}{\bolds\Gamma
}_*\mathbf{Z}^T\mathbf{Z}\bigr)\bigr\Vert_2\bigr\}\\
&&\qquad =(\Delta_1+\Delta_2)\bigl(1+o_p(1)\bigr),
\end{eqnarray*}
where the second inequality holds because ${\mathbf1}_j\mathbf
{N}_n^{-1/2}{\bolds\Gamma}_*=0$ for any $j>q$ by Assumption \ref
{Assum1}, and the
last inequality comes from first-order Taylor expansion and the
definition of $a_n$. In addition, according to the theorem's
condition $\sqrt{n}a_n\stackrel{p}{\rightarrow} 0$, we known
that $\Delta_2$ is $o_p(1)$. Furthermore, based on Lemma~\ref{lema1} and
Assumption~\ref{Assum2}, we have
\begin{eqnarray*}
\sqrt{n}\tr(\mathbf{Z}^T{\mathbf{G}}_n{\bolds\Gamma}_*)&=&\sqrt
{n}\tr(\mathbf{Z}^T{\mathbf{G}}{\bolds\Gamma}_0\mathbf{O})+\sqrt
{n}\tr\bigl(\mathbf{Z}^T({\mathbf{G}}_n\mathbf{N}_n^{1/2}\mathbf
{N}^{-1/2}-{\mathbf{G}}){\bolds{\Gamma
}}_0\mathbf{O}\bigr)\\
&=&\sqrt{n}\tr(\mathbf{Z}^T{\bolds\Gamma}_0\bolds\Lambda_1\mathbf
{O})+\sqrt{n}\tr\bigl(\mathbf{Z}^T({\mathbf{G}}_n-{\mathbf{G}}){\bolds
{\Gamma}}_0\mathbf{O}\bigr)\\
&&{}+\sqrt{n}\tr(\mathbf{Z}^T{\mathbf{G}}_n{\bolds\Gamma}_0\mathbf
{O})\cdot O_p(n^{-1/2})\\
&=&\sqrt{n}\tr\bigl(\mathbf{Z}^T({\mathbf{G}}_n-{\mathbf{G}}){\bolds
{\Gamma}}_0\mathbf{O}\bigr)+O_p(n^{-1/2})\\
&=&\sqrt{n}\tr\bigl(\mathbf{A}^T{\bolds\Gamma}_0^T({\mathbf
{G}}_n-{\mathbf{G}}){\bolds{\Gamma}}_0\mathbf{O}\bigr)\\
&&{}+\sqrt{n}\tr
\bigl(\mathbf{B}^T{\bolds
\Gamma}^T_{0\bot}({\mathbf{G}}_n-{\mathbf{G}}){\bolds{\Gamma
}}_0\mathbf{O}\bigr)+O_p(n^{-1/2})\\
&=&\sqrt{n}\tr\bigl(\mathbf{B}^T{\bolds\Gamma}^T_{0\bot}({\mathbf
{G}}_n-{\mathbf{G}}){\bolds{\Gamma}}_0\bigr)\bigl(1+O_p(n^{-1/2})\bigr),
\end{eqnarray*}
where $\bolds\Lambda=\operatorname{diag}\{\bolds\Lambda_1,\bolds\Lambda
_2\}$ is the
diagonal eigenvalue matrix of $\mathbf G$ with the first $d\times d$
sub-matrix $\bolds\Lambda_1$. By using the definition of $\mathbf
{Z}$ in
(\ref{eqA1}), we get
\begin{eqnarray*}
\tr(\mathbf{Z}^T\mathbf{Z}{\bolds\Gamma}_*^T{\mathbf
{G}}_n{{\bolds
\Gamma}_*})\!-\!\tr(\mathbf{Z}^T\mathbf{G}_n\mathbf{Z})&=&\tr(\mathbf
{Z}^T\mathbf{Z}\mathbf{O}{{\bolds\Gamma}_0^T}{\mathbf{G}}{\bolds
\Gamma}_0\mathbf{O})\!-\!\tr(\mathbf{Z}^T\mathbf{G}\mathbf
{Z})\!+\!O_p(n^{-1/2})\\
&=&\tr(\mathbf{Z}^T\mathbf{Z}\bolds\Lambda_1)-\tr(\mathbf
{Z}^T\mathbf{G}\mathbf{Z})+O_p(n^{-1/2})\\
&=&\tr(\mathbf{A}^T\mathbf{A}\bolds\Lambda_1)+\tr(\mathbf
{B}^T\mathbf{B}\bolds\Lambda_1)-\tr(\mathbf{B}\mathbf{B}^T\bolds
\Lambda_2)\\
&&{}-\tr(\mathbf{A}\mathbf{A}^T\bolds\Lambda_1)+o_p(1)\\
&\geq&(\lambda_d-\lambda_{d+1})\Vert\mathbf{B}\Vert_s^2,
\end{eqnarray*}
where we use the fact $\tr(\mathbf{A}^T\mathbf{A}\bolds\Lambda
_1)-\tr(\mathbf{A}\mathbf{A}^T\bolds\Lambda_1)=0$ because $\mathbf
A$ is skew-symmetric. Here
the last inequality follows from basic properties of trace operator
for semi-positive definite matrix. As a consequence, by the
Cauchy--Schwarz inequality for trace operator, the third term in
$\Delta_1$ is uniformly bounded by $\Vert\mathbf{B}\Vert_s\times
\Vert\sqrt{n}({\mathbf{G}}_n-\mathbf{G}){\bolds{\Gamma}}_0\Vert_s$.
Therefore, as long as the constant $C$ is sufficiently large, the
first two terms in $\Delta_1$ will always dominate the third term
and $\Delta_2$ with arbitrarily large probabilities. This implies inequality~(\ref{eqA2}),
and the proof is complete.
\end{pf*}
\begin{pf*}{Proof of Theorem~\ref{theo2}} (i) To prove this part,
we need represent (\ref{Vtilde}) as vector forms. Define
\begin{eqnarray*}
\mathbf t&=&(\mathbf t_1^T,\ldots,\mathbf t_d^T)^T, \\
h_l(\mathbf t)&=&\mathbf t^T {}\mathbf{C}_l\mathbf t,\qquad l=1,\ldots,d, \\
h_{kl}(\mathbf t)&=&\mathbf t^T\mathbf{C}_{kl}\mathbf t,\qquad (k,l)\in
\mathcal{J},\\
\mathcal{J}&=&\{(k,l)|k,l=1,\ldots,d, k<l\},
\end{eqnarray*}
where $\mathbf t_i$ denotes the $i$th column vector of $\mathbf{V}$,
$\mathbf{C}_l$'s are $pd\times pd$ block-diagonal matrices, $\mathbf
{C}_{kl}$'s
$pd\times pd$ block matrices, $\mathbf{C}_l$ and $\mathbf{C}_{kl}$ contain
$\mathbf{N}_n$ in the $l$th diagonal block and in the $(k,l)$ as well
as $(l,k)$ blocks, respectively. The $pd\times pd$ symmetric
matrices $\mathbf{C}_{kl}$ are defined for all the pairs of different
indices belonging to $\mathcal J$, given by the $d(d-1)/2$ combinations
of the indices $1,\ldots,d$.

By this notation, we have
\[
Q(\bolds\Gamma):=Q^*({\mathbf t})=-\mathbf t^T\mathbf{A}\mathbf
t+\sum_{i=1}^p\theta_i\Vert\mathbf v_i\Vert_2,
\]
where $\mathbf{A}$ is a $pd\times pd$ block-diagonal matrix with all
diagonal blocks $\mathbf{M}_n$. Of course, in the above equation each
$\mathbf v_i$ is regarded as a function of $\mathbf t$.

By using the equality representation of the compact Stiefel
manifolds $\operatorname{St}(p,d)$, (\ref{Vtilde}) is equivalent to
%
%
\begin{eqnarray}\label{eqA3}
&&\min_{\mathbf t}-\Biggl\{\mathbf t^T\mathbf{A}\mathbf
t+\sum_{i=1}^p\theta_i\Vert\mathbf v_i\Vert_2\Biggr\}\qquad\qquad\qquad\nonumber\\[-8pt]\\[-8pt]
&&\eqntext{\mbox{subject to } h_l(\mathbf t)=1, l=1\in[1,d]
\mbox{ and } h_{kl}(\mathbf t)=0, (k,l)\in\mathcal{J}.}
\end{eqnarray}
As a consequence, this enables us to apply an improved global
lagrange multiplier rule proposed by \citet{R97}.

We start by supposing that $\tilde{\mathbf v}_j\neq0$ for all $j$.
According to Theorem 15.2.1 in \citet{R97} [or
Theorem 3.1 in \citet{R02}], a necessary
condition that $\tilde{\mathbf t}_n$ ($\tilde{\mathbf{V}}_n$) is a local
minimum of (\ref{eqA3}) [equation (\ref{Vtilde})] is that, the geodesic gradient
vector of the improved Lagrangian function of (\ref{eqA3}) evaluated at
$\tilde{\mathbf t}_n$ equals to zero. That is,
%
%
\begin{eqnarray}\label{eqA4}
\frac{\partial^g Q^*(\mathbf t)}{\partial\mathbf t}\bigg|_{\mathbf
t=\tilde{\mathbf
t}_n}&\equiv&\biggl[\frac{\partial Q^*(\mathbf t)}{\partial\mathbf
t}-\mathbf{U}(\mathbf{U}'\mathbf{U})^{-1}\mathbf{U}\,\frac{\partial
Q^*(\mathbf
t)}{\partial\mathbf t}\biggr]\bigg|_{\mathbf t=\tilde{\mathbf
t}_n}\nonumber\\[-8pt]\\[-8pt]
:\!&=&\frac{\partial^g f(\mathbf{V}_n)}{\partial\mathbf t}\bigg|_{\mathbf
t=\tilde{\mathbf t}_n}+\frac{\partial^g \rho(\mathbf
{V}_n)}{\partial\mathbf
t}\bigg|_{\mathbf t=\tilde{\mathbf t}_n}=\mathbf0,\nonumber
\end{eqnarray}
where
\[
\mathbf{U}=(\mathbf{C}_1\mathbf t, \ldots,\mathbf{C}_d\mathbf
t,\mathbf{C}_{12}\mathbf t,\mathbf{C}_{13}\mathbf t,\ldots,\mathbf
{C}_{d-1d}\mathbf t)
\]
is a
$(pd\times[d(d+1)/2])$-dimensional matrix, and $\partial^g f(\mathbf
{V}_n)/\partial\mathbf t$ and $\partial^g \rho(\mathbf
{V}_n)/\partial\mathbf t$
are defined in a similar form of $\partial^g Q^*(\mathbf t)/\partial
\mathbf
t$ by replacing $Q^*$ with $f$ and $\rho$, respectively. By Theorem
\ref{theo1} and noting that $\partial f(\mathbf{V}_n)/\partial\mathbf
t$ is
linear in
$\mathbf t$,
\[
\frac{\partial^g f(\mathbf{V}_n)}{\partial\mathbf t}\bigg|_{\mathbf
t=\tilde{\mathbf
t}_n}=\frac{\partial^g f(\mathbf{V}_n)}{\partial\mathbf
t}\bigg|_{\mathbf
t=\widehat{\mathbf t}_n}+O_p(n^{-1/2}),
\]
where $\widehat{\mathbf t}_n$ is the vector form of $\widehat{\mathbf{V}}_n$.
Using Theorem 3.1 of \citet{R02}, we have
$\partial^g f( \mathbf{V}_n)/\partial\mathbf t|_{\mathbf t=\widehat
{\mathbf
t}_n}=\mathbf0$, which yields that $\partial^g f(\mathbf
{V}_n)/\partial\mathbf
t|_{\mathbf t=\tilde{\mathbf t}_n}=O_p(n^{-1/2})$ and as a consequence
\[
\partial^g \rho(\mathbf{V}_n)/{\partial\mathbf t}|_{\mathbf
t=\tilde{\mathbf
t}_n}=O_p(n^{-1/2}).
\]

On the other hand,
\[
\frac{\partial^g \rho(\mathbf{V}_n)}{\partial\mathbf t}\bigg|_{\mathbf
t=\tilde{\mathbf t}_n}=[\mathbf{I}_{pd}-\mathbf{U}(\mathbf
{U}'\mathbf{U})^{-1}\mathbf{U}]\tilde{\bolds\theta}\equiv\mathbf
{H}\tilde{\bolds\theta},
\]
where
\[
\tilde{\bolds\theta}=\biggl(\frac{\theta_1\tilde{t}_{n11}}{\Vert\tilde
{\mathbf
v}_{n1}\Vert_2},\ldots,\frac{\theta_p\tilde{t}_{n1p}}{\Vert\tilde
{\mathbf
v}_{np}\Vert_2},\ldots,\frac{\theta_1\tilde{t}_{nd1}}{\Vert\tilde
{\mathbf
v}_{n1}\Vert_2},\ldots,\frac{\theta_p\tilde{t}_{ndp}}{\Vert\tilde
{\mathbf
v}_{np}\Vert_2}\biggr)^T.
\]
By using the fact that $\mathbf{U}$ has full
column rank and $\mathbf{H}\mathbf{U}=\mathbf0$, we know $\tilde
{\bolds\theta}$
can be expressed through a linear combination of the columns of
$\mathbf{U}$ in probability, that is,
\begin{eqnarray*}
\tilde{\bolds\theta}&=&(\kappa_1\mathbf{C}_1+\cdots+\kappa
_d\mathbf{C}_d
+\kappa_{12}\mathbf{C}_{12}+\kappa_{13}\mathbf{C}_{13}+\cdots+\kappa_{d-1d}\mathbf{C}_{d-1d})
\frac{\Vert\tilde{\bolds
\theta}\Vert_2}{\Vert\tilde{\mathbf
t}_n\Vert_2}\tilde{\mathbf t}_n\\
&&{}+O_p(n^{-1/2}),
\end{eqnarray*}
where $\kappa_1,\ldots,\kappa_{d-1d}$ are a sequence of constants
satisfy they are not all the zeros. Define a sequence of
$pd$-dimensional vectors $\mathbf z_{ij}$'s,
\[
\mathbf z_{ij}=({\mathbf
0}^T,\ldots,\tilde{\mathbf t}_{ni}^T,\ldots,\mathbf{0}^T,\ldots,
\tilde{\mathbf
t}_{ni}^T,\ldots,\mathbf{0}^T)^T,\vadjust{\goodbreak}
\]
for $j \geq i$, say, its
$[(i-1)p+1]$th to the $[(i-1)p+p]$th elements\vspace*{1pt} and $[(j-1)p+1]$th
to the $[(j-1)p+p]$th elements are both $\tilde{\mathbf t}_{ni}$. It is
straightforward to see
%
%
\begin{eqnarray}\label{eqA5}
\kappa_0\kappa_i&=&\mathbf
z_{ii}^T\tilde{\bolds\theta}+O_p(n^{-1/2}),\nonumber\\[-8pt]\\[-8pt]
\kappa_0(\kappa_i+\kappa_{ij})&=&\mathbf
z_{ij}^T\tilde{\bolds\theta}+O_p(n^{-1/2})\qquad \mbox{for }
j>i,\nonumber
\end{eqnarray}
where we denote $\kappa_0=\Vert\tilde{\bolds\theta}\Vert_2/{\Vert
\tilde{\mathbf
t}_n\Vert_2}$. By Theorem~\ref{theo1}, $\tilde{\mathbf
v}_{nj}=O_p(n^{-1/2})$ for
$j>q$. Thus, by recalling the theorem's condition on $a_n$ and
$b_n$, it can be easily verified that (\ref{eqA5}) leads to
\begin{eqnarray*}
\kappa_i+\kappa_{ij}&=&\kappa_0^{-1}\bigl(\mathbf
z_{ij}^T\tilde{\bolds\theta}+O_p(n^{-1/2})\bigr)\\
&\leq& O_p(b_n^{-1})\cdot O_p(a_n + b_n n^{-1/2}+n^{-1/2}) \\
&=& o_p(1).
\end{eqnarray*}
Similarly, $\kappa_i=o_p(1)$. Consequently, we can conclude all the
$\kappa_i$ and $\kappa_{ij}$ equal to zero in probability which
yields contradiction. As a result, with probability tending to 1
(w.p.1), (\ref{eqA4}) cannot hold, which implies there exists $j>q$ so that
\[
\Pr(\tilde{\mathbf v}_{nj}=0)\rightarrow1.
\]

Without loss of generality, we assume $\Pr(\tilde{\mathbf
v}_{np}=0)\rightarrow1$. Let $\mathbf{M}_{n1}$ and $\mathbf{N}_{n1}$ be
the first $(p-1)\times(p-1)$ sub-matrices of $\mathbf{M}_n$ and
$\mathbf{N}_n$, respectively, and $\tilde{\mathbf{V}}_{n1}$ be the
first $p-1$ rows
of $\tilde{\mathbf{V}}_n$. As stated before, $\tilde{\mathbf{V}}_n$
is a local
minimum of the objective function
\[
Q(\mathbf{V};\mathbf{M}_n) = -\tr({\mathbf{V}}^T\mathbf
{M}_n\mathbf{V})+
\sum_{i=1}^{p}\theta_i\Vert{\mathbf v}_i\Vert_2\qquad \mbox{subject to }
{\mathbf{V}}^T\mathbf{N}_n\mathbf{V}=\mathbf{I}_d.
\]
We will show that w.p.1 $\tilde{\mathbf{V}}_{n1}$ is also a local minimum
of the objective function
%
%
\begin{eqnarray}\label{eqA6}
&&Q(\mathbf{V}_{1};\mathbf{M}_{n1})=-\tr({\mathbf{V}_{1}}^T\mathbf
{M}_{n1}\mathbf{V}_{1})+\sum_{i=1}^{p-1}\theta_i\Vert{\mathbf
v}_i\Vert_2\nonumber\\[-8pt]\\[-8pt]
&&\eqntext{\mbox{subject to } {\mathbf{V}_{1}}^T\mathbf{N}_{n1}\mathbf
{V}_{1}=\mathbf{I}_d,}
\end{eqnarray}
w.p.1. Denote the set $\mathcal{A}_{1} =\{ \mathbf{V}_{1} | \Vert
\mathbf{V}_{1} -\tilde{\mathbf{V}}_{n1}\Vert_s < \delta; \mathbf
{V}_{1}^T\mathbf{N}_{n1}\mathbf{V}_{1}=\mathbf{I}_d \}$. For any
$\mathbf{A}_1 \in\mathcal
{A}_{1}$, denote $\mathbf{A} = (\mathbf{A}_1^T, \mathbf{0}^T)^T$. It
is clear
that $\mathbf{A}^T\mathbf{N}_n\mathbf{A}=\mathbf{I}_d$. Given
$\delta$ small
enough, we will have $Q(\mathbf{A};\mathbf{M}_n) \geq Q( \tilde
{\mathbf{V}}_n;\mathbf{M}_n)$ since $\tilde{\mathbf{V}}_n$ is the
local minimum. Note
that $Q(\mathbf{A};\mathbf{M}_n)= Q(\mathbf{A}_1;\mathbf{M}_{n1})$ and
$Q({\tilde{\mathbf{V}}};\mathbf{M}_n)=Q(\tilde{\mathbf
{V}}_{n1};\mathbf{M}_{n1})$
w.p.1. Consequently, we have
\[
Q(\mathbf{A}_1;\mathbf{M}_{n1}) \geq Q( \tilde{\mathbf
{V}}_{n1};\mathbf{M}_{n1})\qquad \mbox{w.p.1},
\]
for all $\mathbf{A}_1 \in\mathcal{A}$
provided that $\delta$ is sufficiently small. Hence, we can conclude
that $\tilde{\mathbf{V}}_{n1}$ is also a local minimum of the objective
function $Q(\mathbf{V}_{1};\mathbf{M}_{n1})$ w.p.1.

Rewriting (\ref{eqA6}) as a similar form to (\ref{eqA3}) and following the same
arguments above in proving $\Pr(\tilde{\mathbf v}_{np}=0)\rightarrow1$,
we can show that there exists $q<j<p$ so that $\Pr(\tilde{\mathbf
v}_{nj}=0)\rightarrow1$. The remaining proofs can be completed by
deduction.

(ii) For convenience purposes, first decompose the matrix
$\mathbf{M}_{n}$ and $\mathbf{N}_{n}$ into the following block form:
\begin{eqnarray*}
\mathbf{M}_n=\left[
\matrix{\mathbf{M}_{n(q)}& \mathbf{M}_{12}\cr\mathbf{M}_{21}
&\mathbf{M}_{n(p-q)}}
\right],\qquad \mathbf{N}_n=\left[
\matrix{\mathbf{N}_{n(q)}& \mathbf{N}_{12}\cr\mathbf{N}_{21}
&\mathbf{N}_{n(p-q)}}
\right],
\end{eqnarray*}
where $\mathbf{M}_{n(q)}$ and $\mathbf{N}_{n(q)}$ are the first
$q\times
q$ sub-matrices. It then follows that
\[
f(\mathbf{V};\mathbf{M}_n)=-\tr\bigl(\mathbf{V}^T_{(q)}\mathbf
{M}_{n(q)}\mathbf{V}_{(q)}\bigr)-\tr\bigl(\mathbf{V}^T_{(p-q)}\mathbf
{M}_{n(p-q)}\mathbf{V}_{(p-q)}\bigr).
\]
Next we will show ${\tilde{\mathbf{V}}_{n(q)}}=\widehat{\mathbf
{V}}_{n(O)}(1+o_p(n^{-1/2}))$. Similar to the proof of Theorem~\ref{theo1},
since ${\tilde{\mathbf{V}}_{n(p-q)}}=0$ w.p.1, it suffices to show, for
any arbitrarily small $\varepsilon>0$, there exits a sufficiently
large constant $C$, such that
%
%
\begin{eqnarray}\label{eqA7}
&&\lim_{n}\inf\pr\Bigl(\inf_{\mathbf{Z}\in
T_{\widehat{\bolds\Gamma}_{n(O)}}(q,d)\dvtx\Vert\mathbf{B}\Vert
_s=C}Q\bigl(R\bigl({\widehat{\bolds\Gamma}_{n(O)}}+a_n\mathbf{Z}\bigr);\mathbf
{G}_{n(q)},\mathbf{N}_{n(q)}\bigr)\nonumber\\
&&\qquad\quad\hspace*{136.1pt}{}
>Q\bigl(\widehat{\bolds\Gamma}_{n(O)};\mathbf{G}_{n(q)},\mathbf
{N}_{n(q)}\bigr)\Bigr)\\
&&\qquad>1-\varepsilon,\nonumber
\end{eqnarray}
where
\[
{\widehat{\bolds\Gamma}_{n(O)}}=\mathop{\arg\min}_{\bolds\Gamma
\in
\mathbb{R}^{q\times d}}-\tr\bigl(\bolds\Gamma^T\mathbf{G}_{n(q)}\bolds
\Gamma\bigr)\qquad
\mbox{subject to } \bolds\Gamma^T\bolds\Gamma=\mathbf{I}_d
\]
and $\mathbf{G}_{n(q)}=\mathbf{N}_{n(q)}^{-1/2}\mathbf
{M}_{n(q)}\mathbf{N}_{n(q)}^{-1/2}$. Note that
\begin{eqnarray*}
&&
a_n^{-2}\bigl\{Q\bigl(R\bigl(\widehat{\bolds\Gamma}_{n(O)}+a_n\mathbf{Z}\bigr);\mathbf
{G}_{n(q)},\mathbf{N}_{n(q)}\bigr)-Q\bigl(\widehat{\bolds\Gamma
}_{n(O)};\mathbf{G}_{n(q)},\mathbf{N}_{n(q)}\bigr)\bigr\}\\
&&\qquad\geq\bigl[-\tr\bigl(\mathbf{Z}^T\mathbf{G}_{n(q)}\mathbf{Z}\bigr)-2a_n^{-1}\tr
\bigl(\mathbf{Z}^T\mathbf{G}_{n(q)}{\widehat{\bolds\Gamma}}_{n(O)}\bigr)+\tr
\bigl(\mathbf{Z}^T\mathbf{Z}{\widehat{\bolds\Gamma}}{}^T_{n(O)}{\mathbf
{G}}_{n(q)}{{\widehat{\bolds\Gamma
}}_{n(O)}}\bigr)\bigr]\\
&&\qquad\quad{}\times\bigl(1+o_p(1)\bigr)\\
&&\qquad\quad{}-q\bigl\Vert{\mathbf1}_j\mathbf{N}_{n(q)}^{-1/2}\bigl(\mathbf
{Z}-(1/2)a_n\widehat{\bolds\Gamma}_{n(O)}\mathbf{Z}^T\mathbf
{Z}\bigr)\bigr\Vert_2,
\end{eqnarray*}
where $2a_n^{-1}\tr(\mathbf{Z}^T\mathbf{G}_{n(q)}{\widehat{\bolds
\Gamma}}_{n(O)})=0$ by using Lemma~\ref{lema2}, and
\[
-\tr\bigl(\mathbf{Z}^T\mathbf{G}_{n(q)}\mathbf{Z}\bigr)+\tr\bigl(\mathbf
{Z}^T\mathbf{Z}{\widehat{\bolds\Gamma}}{}^T_{n(O)}{\mathbf
{G}}_{n(q)}{{\widehat{\bolds\Gamma
}}_{n(O)}}\bigr)>0.
\]
Using the similar arguments in the proof of Theorem~\ref{theo1}, we can show
(\ref{eqA7}) holds. This\vadjust{\goodbreak} implies that $\sqrt{n}\tilde{\bolds\Gamma
}_{n(q)}$ is
asymptotically equivalent to $\sqrt{n}\widehat{\bolds\Gamma}_{n(O)}$
where
\[
{\tilde{\bolds\Gamma}_{n(q)}}=\mathop{\arg\min}_{\bolds\Gamma
\in
\mathbb{R}^{q\times d}}Q\bigl(\bolds\Gamma;\mathbf{G}_{n(q)},\mathbf{N}_{n(q)}\bigr)
\qquad\mbox{subject to } \bolds\Gamma^T\bolds\Gamma=\mathbf{I}_d,
\]
and thus it
follows that $\sqrt{n}D(\mathbf{N}_{n(q)}^{1/2}\tilde{\mathbf
{V}}_{n(q)},\mathbf{N}_{n(q)}^{1/2}\widehat{\mathbf
{V}}_{n(O)})=o_p(1)$ which
completes the proof.
\end{pf*}
\begin{pf*}{Proof of Proposition~\ref{prop3}}
To illustrate the idea, we
elaborate on verifying the condition (\ref{oc}) for DR. In this
case, by equation (5) in \citet{LiWang2007}, $\mathbf{M}_n$ can be
reexpressed as
\begin{eqnarray*}
\mathbf{M}_n&=&2\bigl\{{\bolds\Sigma}_n^{1/2}\widehat{E}[\widehat{\var
}(\mathbf{z}|\tilde{y})-\mathbf{I}_p]^2{\bolds\Sigma
}_n^{1/2}\\
&&\hspace*{8.5pt}{}+{\bolds\Sigma
}_n^{1/2}\widehat{E}\bigl[\bigl(\widehat{\var}(\mathbf{z}|\tilde{y})-\mathbf
{I}_p\bigr)\widehat{E}(\mathbf{z}|\tilde{y})\widehat{E}(\mathbf{z}^T
|\tilde
{y})\bigr]{\bolds\Sigma}_n^{1/2}\\
&&\hspace*{8.5pt}{}+{\bolds\Sigma}^{1/2}_n\widehat{E}\bigl[\widehat{E}(\mathbf{z}|\tilde
{y})\widehat{E}(\mathbf{z}^T |\tilde{y})\bigl(\widehat{\var}(\mathbf
{z}|\tilde{y})-\mathbf{I}_p\bigr)\bigr]{\bolds\Sigma}_n^{1/2}\\
&&\hspace*{8.5pt}{}+{\bolds
\Sigma}_n^{1/2}\widehat{E}[\widehat{E}(\mathbf{z}|\tilde
{y})\widehat{E}(\mathbf{z}^T |\tilde{y})]^2{\bolds
\Sigma}_n^{1/2}\\
&&\hspace*{8.5pt}{} + {\bolds
\Sigma}_n^{1/2}\widehat{E}^2[\widehat{E}(\mathbf{z}|\tilde
{y})\widehat{E}(\mathbf{z}^T |\tilde{y})]{\bolds
\Sigma}_n^{1/2}\\
&&\hspace*{8.5pt}{}+ {\bolds
\Sigma}_n^{1/2}\widehat{E}[\widehat{E}(\mathbf{z}^T|\tilde
{y})\widehat{E}(\mathbf{z}|\tilde{y})]\widehat{E}[\widehat
{E}(\mathbf{z}|\tilde{y})\widehat{E}(\mathbf{z}^T |\tilde
{y})]{\bolds
\Sigma}_n^{1/2}\bigr\}\\
:\!&=&2(\mathbf{M}_{n1}+\cdots+\mathbf{M}_{n6}).
\end{eqnarray*}
Here, $\tilde{y}$ is the discretized $y$ over a collection of slices,
$\widehat{\var}(\mathbf{z}|\tilde{y})$ denotes the sample covariance
matrix of $\mathbf{z}$ within a slice, $\widehat{E}(\cdot)$ denotes the
weighted average across slices. Next, we will show $\mathbf
{M}_{n(O)i}=\mathbf{M}_{n(q)i}+O_p(n^{-1})$ for $i=1,\ldots,6$.

Now we first deal with $\mathbf{M}_{n1}$. Rewrite it as
\[
\mathbf{M}_{n1}=\widehat{E}\{[\widehat{\var}({\mathbf x}|\tilde
{y})-\bolds\Sigma_n]\bolds\Sigma_n^{-1}[\widehat{\var}({\mathbf
x}|\tilde{y})-\bolds\Sigma_n]\}.
\]
We assume that the collection of slices is fixed; that is, it does
not vary with $n$. This implies that the sample conditional moments
such as $\widehat{\var}({\mathbf x}|\tilde{y})$ are
$\sqrt{n}$-consistent estimates of their population-level
counterparts, such as ${\var({\mathbf x}|\tilde{y})}$. Let $\bolds
\Omega$
be the matrix consisting of the first $q$ columns of the matrix
$\mathbf{I}_p$. Then, by definition,
\begin{eqnarray*}
\mathbf{M}_{n(O)1}&=&\bolds\Omega^{T}\widehat{E}\{[\widehat{\var
}({\mathbf x}|\tilde{y})-\bolds\Sigma_n]\bolds\Omega(\bolds\Omega
^{T}\bolds\Sigma_n\bolds\Omega
)^{-1}\bolds\Omega^{T}[\widehat{\var}({\mathbf x}|\tilde
{y})-\bolds\Sigma_n]\}\bolds\Omega,\\
\mathbf{M}_{n(q)1}&=&\bolds\Omega^{T}\widehat{E}\{[\widehat{\var
}({\mathbf x}|\tilde{y})-\bolds\Sigma_n]\bolds\Sigma
_n^{-1}[\widehat{\var}({\mathbf x}|\tilde{y})-\bolds\Sigma_n]\}
\bolds\Omega.
\end{eqnarray*}
Let ${\mathbf
P}_{\bolds\Omega}(\bolds\Sigma_n)=\bolds\Omega(\bolds\Omega
^{T}\bolds\Sigma_n\bolds\Omega
)^{-1}\bolds\Omega^{T}\bolds\Sigma_n$
and let ${\mathbf Q}_{\bolds\Omega}(\bolds\Sigma_n)=\mathbf
{I}_p-{\mathbf
P}_{\bolds\Omega}(\bolds\Sigma_n)$. Then
\begin{eqnarray*}
\mathbf{M}_{n(q)1}&=&\bolds\Omega^{T}\widehat{E}\{[\widehat{\var
}({\mathbf x}|\tilde{y})-\bolds\Sigma_n][{\mathbf
P}_{\bolds\Omega}(\bolds\Sigma_n)+\mathbf{
Q}_{\bolds\Omega}(\bolds\Sigma_n)]\bolds\Sigma_n^{-1}\\
&&\hspace*{27.2pt}{}\times[\mathbf{P}_{\bolds\Omega}(\bolds\Sigma_n)+\mathbf{Q}_{\bolds
\Omega}(\bolds\Sigma_n)]^T[\widehat{\var}({\mathbf x}|\tilde
{y})-\bolds\Sigma_n]\}\bolds\Omega\\
:\!&=&\widehat{E}(\mathbf
{M}_{1\mathrm{I}}+\mathbf{M}_{1\mathrm{II}}+\mathbf{M}_{1\mathrm{III}}+\mathbf{M}_{1\mathrm{IV}}),
\end{eqnarray*}
where
\begin{eqnarray*}
\mathbf{M}_{1\mathrm{I}}&=&\bolds\Omega^{T}[\widehat{\var}({\mathbf
x}|\tilde{y})-\bolds\Sigma_n]\mathbf{P}_{\bolds\Omega}(\bolds
\Sigma_n)\bolds\Sigma_n^{-1}\mathbf{P}_{\bolds\Omega}^T(\bolds
\Sigma_n)[\widehat{\var}({\mathbf x}|\tilde{y})-\bolds\Sigma
_n]\bolds\Omega,\\
\mathbf{M}_{1\mathrm{II}}&=&\bolds\Omega^{T}[\widehat{\var}({\mathbf
x}|\tilde{y})-\bolds\Sigma_n]\mathbf{Q}_{\bolds\Omega}(\bolds
\Sigma_n)\bolds\Sigma_n^{-1}\mathbf{P}_{\bolds\Omega}^T(\bolds
\Sigma_n)[\widehat{\var}({\mathbf x}|\tilde{y})-\bolds\Sigma
_n]\bolds\Omega,\\
\mathbf{M}_{1\mathrm{III}}&=&\bolds\Omega^{T}[\widehat{\var}({\mathbf
x}|\tilde{y})-\bolds\Sigma_n]\mathbf{P}_{\bolds\Omega}(\bolds
\Sigma_n)\bolds\Sigma_n^{-1}\mathbf{Q}_{\bolds\Omega}^T(\bolds
\Sigma_n)[\widehat{\var}({\mathbf x}|\tilde{y})-\bolds\Sigma
_n]\bolds\Omega,\\
\mathbf{M}_{1\mathrm{IV}}&=&\bolds\Omega^{T}[\widehat{\var}({\mathbf
x}|\tilde{y})-\bolds\Sigma_n]\mathbf{Q}_{\bolds\Omega}(\bolds
\Sigma_n)\bolds\Sigma_n^{-1}\mathbf{Q}_{\bolds\Omega}^T(\bolds
\Sigma_n)[\widehat{\var}({\mathbf x}|\tilde{y})-\bolds\Sigma
_n]\bolds\Omega.
\end{eqnarray*}
It can be easily seen that $\widehat{E}(\mathbf{M}_{1\mathrm{I}})$ is exactly
$\mathbf{M}_{n(O)1}$. We will show that $\mathbf{M}_{1\mathrm{II}}$, $\mathbf
{M}_{1\mathrm{III}}$ and $\mathbf{M}_{1\mathrm{IV}}$ are of the order $O_p(n^{-1})$.
Note that
\begin{eqnarray*}
&&\mathbf{Q}_{\bolds\Omega}^T(\bolds\Sigma_n)[\widehat{\var
}({\mathbf x}|\tilde{y})-\bolds\Sigma_n]\\
&&\qquad=[\mathbf{Q}_{\bolds\Omega}^T(\bolds\Sigma
)+O_p(n^{-1/2})][\widehat{\var}({\mathbf x}|\tilde{y})-\bolds\Sigma
+O_p(n^{-1/2})]\\
&&\qquad=\mathbf{Q}_{\bolds\Omega}^T(\bolds\Sigma)[\widehat{\var
}({\mathbf x}|\tilde{y})-\bolds\Sigma]+O_p(n^{-1/2}).
\end{eqnarray*}
By construction, $\mathcal{S}_{y|{\mathbf x}}\subseteq
\operatorname{span}(\bolds\Omega)$. Under certain conditions
[\citet{Cook98a}],
we know $\operatorname{span}\{\bolds\Sigma^{-1}[\bolds\Sigma
-\var
({\mathbf x}|y)]\}\subseteq\mathcal{S}_{y|{\mathbf x}}$. Hence,
\[
\operatorname{span}\{\bolds\Sigma^{-1}[\bolds\Sigma-\var({\mathbf
x}|y)]\}
\subseteq\operatorname{span}(\bolds\Omega).
\]
It then follows that
%
%
\begin{equation}\label{eqA8}
\mathbf{Q}_{\bolds\Omega}(\bolds\Sigma)\bolds\Sigma^{-1}[{\var
}({\mathbf x}|\tilde{y})-\bolds\Sigma]=\bolds\Sigma^{-1}\mathbf
{Q}_{\bolds\Omega}^T(\bolds\Sigma)[{\var}({\mathbf x}|\tilde
{y})-\bolds\Sigma]=\mathbf{0}.
\end{equation}
Thus, we have $\mathbf{M}_{1\mathrm{IV}}=O_p(n^{-1/2})\cdot
O_p(n^{-1/2})=O_p(n^{-1})$.

Substituting $\mathbf{P}_{\bolds\Omega}(\bolds\Sigma_n)=\mathbf
{I}_p-\mathbf{Q}_{\bolds\Omega}(\bolds\Sigma_n)$ into $\mathbf
{M}_{1\mathrm{II}}$ and using $\mathbf{Q}_{\bolds\Omega}(\bolds\Sigma_n)$'s
idempotency, we have
\[
\mathbf{M}_{1\mathrm{II}}=\bolds\Omega^{T}[\widehat{\var}({\mathbf
x}|\tilde{y})-\bolds\Sigma_n]\mathbf{Q}_{\bolds\Omega}(\bolds
\Sigma_n)\mathbf{Q}_{\bolds\Omega}(\bolds\Sigma_n)\bolds\Sigma
_n^{-1}[\widehat{\var}({\mathbf x}|\tilde{y})-\bolds\Sigma
_n]\bolds\Omega-\mathbf{M}_{1\mathrm{IV}}.
\]
By using (\ref{eqA8}) again, we know that $\mathbf{M}_{1\mathrm{II}}=O_p(n^{-1})$.
Similarly, $\mathbf{M}_{1\mathrm{III}}=O_p(n^{-1})$. From these, we deduce that
$\mathbf{M}_{1\mathrm{II}}$, $\mathbf{M}_{1\mathrm{III}}$, $\mathbf{M}_{1\mathrm{IV}}$ are all of
order $O_p(n^{-1})$. Since $\widehat{E}(\mathbf{M}_{1\mathrm{II}}+\mathbf
{M}_{1\mathrm{III}}+\mathbf{M}_{1\mathrm{IV}})$ is the sum of finite number of terms each
of the order $O_p(n^{-1})$, it is itself of this order. It follows
that $\mathbf{M}_{n(O)1}=\mathbf{M}_{n(q)1}+O_p(n^{-1})$.

Next, let us deal with $\mathbf{M}_{n2}$. Similar to $\mathbf{M}_{n(q)1}$,
$\mathbf{M}_{n(q)2}$ can be divided into four terms $\mathbf
{M}_{n(q)2}=\mathbf{M}_{n(O)2}+\mathbf{M}_{2\mathrm{II}}+\mathbf
{M}_{2\mathrm{III}}+\mathbf{M}_{2\mathrm{IV}}$, where
\begin{eqnarray*}
\mathbf{M}_{2\mathrm{II}}&=&\bolds\Omega^{T}[\widehat{\var}({\mathbf
x}|\tilde{y})-\bolds\Sigma_n]\mathbf{Q}_{\bolds\Omega}(\bolds
\Sigma_n)\bolds\Sigma_n^{-1}\\
&&\hspace*{0pt}{}\times\mathbf{P}_{\bolds\Omega}^T(\bolds
\Sigma_n)\{[\widehat{E}({\mathbf x}|\tilde{y})-\widehat{E}({\mathbf
x})][\widehat{E}({\mathbf x}^T |\tilde
{y})-\widehat{E}({\mathbf x}^T)]\}\bolds\Omega,\\
\mathbf{M}_{2\mathrm{III}}&=&\bolds\Omega^{T}[\widehat{\var}({\mathbf
x}|\tilde{y})-\bolds\Sigma_n]\mathbf{P}_{\bolds\Omega}(\bolds
\Sigma_n)\bolds\Sigma_n^{-1}\\
&&{}\times\mathbf{Q}_{\bolds\Omega}^T(\bolds
\Sigma_n)\{[\widehat{E}({\mathbf x}|\tilde{y})-\widehat{E}({\mathbf
x})][\widehat{E}({\mathbf x}^T |\tilde
{y})-\widehat{E}({\mathbf x}^T)]\}\bolds\Omega,\\
\mathbf{M}_{2\mathrm{IV}}&=&\bolds\Omega^{T}[\widehat{\var}({\mathbf
x}|\tilde{y})-\bolds\Sigma_n]\mathbf{Q}_{\bolds\Omega}(\bolds
\Sigma_n)\bolds\Sigma_n^{-1}\\
&&{}\times\mathbf{Q}_{\bolds\Omega}^T(\bolds
\Sigma_n)\{[\widehat{E}({\mathbf x}|\tilde{y})-\widehat{E}({\mathbf
x})][\widehat{E}({\mathbf x}^T
|\tilde{y})-\widehat{E}({\mathbf x}^T)]\}\bolds\Omega.
\end{eqnarray*}
Under the linearity condition, we know $\operatorname{span}\{
[{E}({\mathbf
x}|\tilde{y})-{E}({\mathbf x})]\}\subseteq\mathcal{S}_{y|{\mathbf
x}}$ [\citet{Cook98a}]. Hence,
\[
\operatorname{span}\{[{E}({\mathbf x}|\tilde{y})-{E}({\mathbf x})]\}
\subseteq
\operatorname{span}(\bolds\Omega).
\]
It then follows that
%
%
\begin{equation}\label{eqA9}
\mathbf{Q}_{\bolds\Omega}(\bolds\Sigma)\bolds\Sigma
^{-1}[{E}({\mathbf x}|\tilde{y})-{E}({\mathbf x})]=\bolds\Sigma
^{-1}\mathbf{Q}_{\bolds\Omega}^T(\bolds\Sigma)[{E}({\mathbf
x}|\tilde{y})-{E}({\mathbf x})]=\mathbf{0}.
\end{equation}
By using (\ref{eqA9}) and the similar arguments for $\mathbf{M}_{n(q)1}$, we
can show that $\mathbf{M}_{2\mathrm{II}}$, $\mathbf{M}_{2\mathrm{III}}$ and $\mathbf{M}_{2\mathrm{IV}}$
are all of order $O_p(n^{-1})$. Thus, we can conclude that $\mathbf
{M}_{n(O)2}=\mathbf{M}_{n(q)2}+O_p(n^{-1})$.

By (\ref{eqA8}) and (\ref{eqA9}), $\mathbf{M}_{n(O)i}=\mathbf{M}_{n(q)i}+O_p(n^{-1})$
for $i=3,\ldots,6$, can be proved in a similar fashion to the
foregoing proofs. We omit the details here for saving some space. It
follows that for the DR method,
\[
\mathbf{M}_{n(O)}=\mathbf{M}_{n(q)}+O_p(n^{-1}).
\]
Thus,   condition
(\ref{oc}) is satisfied as long as $(na_n)^{-1}=O_p(1)$.

Note that for SAVE, $\mathbf{M}_n$ takes the form of $\mathbf
{M}_{n1}$ for
DR. Thus,   condition~(\ref{oc}) is also satisfied for SAVE.
\end{pf*}
\end{appendix}

\section*{Acknowledgments}

The authors thank the Associate Editor and two anonymous
referees for their many helpful comments that have resulted in
significant improvements in the article. Specially, we are grateful
to the Associate Editor for helping us complete the proof of
Proposition~\ref{prop3}. The
authors would also like to thank Dr. Jianhui Zhou and Dr. Liqiang Ni
for providing us the codes for computing the $C^3$ and SSIR
estimators.


%
\printaddresses

\end{document}